\documentclass[12pt]{article}

\usepackage{subcaption}
\usepackage{bigints}

\usepackage{float}

\usepackage{multirow}

\usepackage{tikz}

\usepackage{algorithm}
\usepackage{algpseudocode}

\algrenewcommand\algorithmicrequire{\textbf{Input:}}
\algrenewcommand\algorithmicensure{\textbf{Output:}}

\usetikzlibrary{calc,backgrounds,arrows,matrix}

\usepackage{enumitem}

\usepackage{listings}

\usepackage{color}
\definecolor{lightblue}{rgb}{0,0.2,0.5}

\usepackage[colorlinks=true, urlcolor=lightblue,linkcolor=lightblue, citecolor=lightblue]{hyperref}

\usepackage{amssymb,amsmath}

\usepackage{color}
\usepackage{mathrsfs}

\usepackage{mathtools}
\usepackage{accents}

\DeclareMathAlphabet{\eufrak}{U}{}{}{}

\SetMathAlphabet\eufrak{normal}{U}{euf}{m}{n}
\SetMathAlphabet\eufrak{bold}{U}{euf}{b}{n}

\DeclareMathOperator*{\argmin}{arg\,min}

\allowdisplaybreaks

 \def\qu{{\mathord{\mathbb Z}}}

 \def\Var{{\mathrm{{\rm Var}}}}

 \def\sZZ{{\rm Z\kern-.45em{}Z}}

 \def\sQQ{{\kern 0.27em \vrule height1.45ex width0.03em depth0em
           \kern-0.30em \rm Q}}
 \def\qu{{\mathchoice
         {\sQQ}
         {\sQQ}
   {\kern 0.225em \vrule height1.05ex width0.025em depth0em \kern-0.25em \rm Q}
   {\kern 0.180em \vrule height0.78ex width0.020em depth0em \kern-0.20em \rm Q}
         }}
 \def\sGG{{\kern 0.27em \vrule height1.45ex width0.03em depth0em
           \kern-0.30em \rm G}}
 \def\gg{{\mathchoice
         {\sGG}
         {\sGG}
   {\kern 0.225em \vrule height1.05ex width0.025em depth0em \kern-0.25em \rm G}
   {\kern 0.180em \vrule height0.78ex width0.020em depth0em \kern-0.20em \rm G}
         }}

 \newtheorem{prop}{Proposition}[section]
 
 \newtheorem{definition}[prop]{Definition}

 \newtheorem{remark}[prop]{Remark}

\newcommand{\abs}[1]{\lvert#1\rvert}

\newcommand{\re}{\mathrm{e}}

 \newcounter{hyp}
\newenvironment{Proofy}{\removelastskip\par\medskip
\noindent{\em Proof of Proposition} \rm}{\penalty-20\null\hfill$\square$\par\medbreak}

\newenvironment{Proof}{\removelastskip\par\medskip \noindent{\em Proof.} \rm}{\penalty-20\null\hfill$\square$\par\medbreak}

\def\bprf{\begin{Proof}}
\def\nprf{\end{Proof}}
\def\bdes{\begin{description}}
\def\ndes{\end{description}}

\newtheorem{thm}{Theorem}[section]

\def\bdef{\begin{defn}}
\def\ndef{\end{defn}}
\def\bthm{\begin{thm}}
\def\nthm{\end{thm}}
\def\bprop{\begin{prop}}
\def\nprop{\end{prop}}
\def\brmk{\begin{remark}}
\def\nrmk{\end{remark}}
\def\bexa{\begin{exa}}
\def\nexa{\end{exa}}
\def\blem{\begin{lem}}
\def\nlem{\end{lem}}
\def\bcor{\begin{cor}}
\def\ncor{\end{cor}}
\def\bexe{\begin{exe}}
\def\nexe{\end{exe}}

\newcommand{\E}{\mathbb{E}}

\newcommand{\nn}{\mathbb{N}}
\newcommand{\real}{\mathbb{R}}

\newcommand{\zz}{\mathbb{Z}}

\def\Var{\mathop{\hbox{\rm Var}}\nolimits}

\def\og{\leavevmode\raise.3ex
     \hbox{$\scriptscriptstyle\langle\!\langle$~}}
\def\fg{\leavevmode\raise.3ex
     \hbox{~$\!\scriptscriptstyle\,\rangle\!\rangle$}~}

\title{\Huge
 Numerical solution of the incompressible Navier-Stokes equation by a deep branching algorithm
}

\author{
 Jiang Yu Nguwi\footnote{\href{mailto:nguw0003@e.ntu.edu.sg}{nguw0003@e.ntu.edu.sg}
 }
 \qquad
 Guillaume Penent\footnote{\href{mailto:PENE0001@e.ntu.edu.sg}{pene0001@e.ntu.edu.sg}}
 \qquad Nicolas Privault\footnote{
\href{mailto:nprivault@ntu.edu.sg}{nprivault@ntu.edu.sg}
 }
 \\
  \small
Division of Mathematical Sciences
\\
\small
School of Physical and Mathematical Sciences
\\
\small
Nanyang Technological University
\\
\small
21 Nanyang Link, Singapore 637371
}

\allowdisplaybreaks

\usepackage{cases}

\usepackage{empheq}
\usepackage{bm}

\usepackage{bbm}

\makeatletter
\newcommand*\rel@kern[1]{\kern#1\dimexpr\macc@kerna}
\newcommand*\widebar[1]{
  \begingroup
  \def\mathaccent##1##2{
    \rel@kern{0.8}
    \overline{\rel@kern{-0.8}\macc@nucleus\rel@kern{0.2}}
    \rel@kern{-0.2}
  }
  \macc@depth\@ne
  \let\math@bgroup\@empty \let\math@egroup\macc@set@skewchar
  \mathsurround\z@ \frozen@everymath{\mathgroup\macc@group\relax}
  \macc@set@skewchar\relax
  \let\mathaccentV\macc@nested@a
  \macc@nested@a\relax111{#1}
  \endgroup
}
\makeatother

\numberwithin{equation}{section}

\begin{document}

\maketitle

\baselineskip0.6cm

\vspace{-0.6cm}

\begin{abstract}
We present an algorithm for the numerical solution of systems of fully nonlinear PDEs using stochastic coded branching trees. This approach covers functional nonlinearities involving gradient terms of arbitrary orders and it requires only a boundary condition over space at a given terminal time $T$ instead of Dirichlet or Neumann boundary conditions at all times as in standard solvers. Its implementation relies on Monte Carlo estimation, and uses neural networks that perform a meshfree functional estimation on a space-time domain. The algorithm is applied to the numerical solution of the Navier-Stokes equation and is benchmarked to other implementations in the cases of the Taylor-Green vortex and Arnold-Beltrami-Childress flow.
\end{abstract}



\baselineskip0.7cm

\section{Introduction}
This paper is concerned with the numerical solution of
systems of $d+1$ fully nonlinear
coupled parabolic partial differential equations (PDEs)
and Poisson equations on $[0,T] \times \real^d$, of the form
\begin{equation}
\label{eq:main pde}
\begin{cases}
  \displaystyle
  \partial_t u_i(t,x) + \nu \Delta u_i(t,x)
    + f_i\big(
    \partial_{\bar{\alpha}^1}u_0 (t,x)
    ,
    \ldots ,
    \partial_{\bar{\alpha}^q }u_0 (t,x)
    ,
        \partial_{\bar{\alpha}^{q+1}}u_{\beta_{q+1}}(t,x)
        ,
        \ldots ,
                   \partial_{\bar{\alpha}^n}u_{\beta_n}(t,x)
        \big) = 0,
  \medskip
  \\
  \Delta u_0(t, x)
  = f_0\big(
          \partial_{\bar{\alpha}^{q+1}}u_{\beta_{q+1}}(t,x)
        ,
        \ldots ,
                   \partial_{\bar{\alpha}^n}u_{\beta_n}(t,x)
        \big),
  \medskip
  \\
u_i(T,x) = \phi_i (x),
\quad (t,x) = (t,x_1, \ldots, x_d) \in [0,T] \times \real^d,
\quad
 i = 1,\ldots , d,
\end{cases}
\end{equation}
where $q\in \{0,1,\ldots , n\}$,
$\partial_t u (t,x) = \partial u(t,x) / \partial t$,
$\nu > 0$, $\Delta  = \sum\limits_{i = 1}^d \partial^2 / \partial x_i^2$
            is the standard $d$-dimensional Laplacian,
                          $1 \leq \beta_j \leq d$ for $q<j\leq n$,
$\bar{\alpha}^i = (\alpha^i_1, \dots, \alpha^i_d ) \in \nn^d$,
$i = 0, 1, \ldots, n$, $f_i(x_1,\ldots , x_n)$ and
$f_0(x_{q+1},\ldots , x_n)$
 are smooth functions of the derivatives
$$
\partial_{\bar{\alpha}^i} u(t,x)
=
\frac{\partial^{\alpha_1^i}}{\partial x_1^{\alpha_1^i}}
\cdots
\frac{\partial^{\alpha_d^i}}{\partial x_d^{\alpha_d^i}}
u(t,x_1,\ldots , x_d), \qquad (x_1,\ldots , x_d)\in \real^d.
$$
We note that the problem~\eqref{eq:main pde} is posed using
the terminal time boundary condition $u_i(T,x) = \phi_i (x)$
in $(x_1, \ldots, x_d) \in \times \real^d$, instead of assuming
Dirichlet and Neumann boundary conditions at all times as is usually
 the case in the finite difference and mesh-based literature. 

\medskip

As is well known, standard numerical schemes for solving \eqref{eq:main pde}
by e.g. finite differences or finite elements
suffer from a high computational cost which typically grows exponentially with the dimension $d$.
This motivates the study of probabilistic representations of \eqref{eq:main pde},
which, combined with meshfree Monte Carlo approximation, 
can overcome the curse of dimensionality. 
In addition, it is not clear how the standard numerical schemes
can be applied when boundary conditions are not available.

\medskip

Probabilistic representations for the solutions of
first and second order nonlinear PDEs 
can be obtained by writing $u(t,x) \in \real$ as
 $u(t,x) = Y_t^{t,x}$,
 where $(Y_s^{t,x})_{t \leq s \leq T}$ is the
 solution of first or second order backward stochastic differential equation (BSDE),
 see \cite{pardouxpeng}, \cite{touzi}, \cite{soner},
 and \cite{han2018solving} for a deep learning implementation.

\medskip

On the other hand, stochastic branching diffusion mechanisms
 (\cite{skorohodbranching}, \cite{inw}, \cite{hpmckean})
have been applied
to the probabilistic representation of the solutions of nonlinear PDEs,
see e.g.\ 
 \cite{henry-labordere2012}, \cite{labordere}
 for the case of polynomial first order gradient nonlinearities,
 and \cite{fahim}, \cite{tanxiaolu}, \cite{guowenjie}, 
\cite{huangshuo} for finite difference schemes combined with Monte Carlo estimation
for fully nonlinear PDEs with gradients of order up to two.
However,
extending the above approaches to nonlinearities in higher order derivatives
involves technical difficulties
linked to the integrability of the Malliavin-type weights used in
repeated integration by parts argument, see page~199 of \cite{labordere}.

\medskip

In \cite{penent2022fully}, a stochastic branching method
that carries information on (possibly functional) nonlinearities along a random tree
has been introduced, with the aim of providing
 Monte Carlo schemes for the numerical solution of
fully nonlinear PDEs with gradients of arbitrary orders on the real line.
This method has been implemented on $\real^d$ in \cite{nguwipenentprivault}
using a neural network approach
to efficiently approximate the PDE solution $u(t, x) \in \real$
over a bounded domain in $[0,T] \times \real^d$.

\medskip

In this paper, we extend the approaches in
 \cite{nguwipenentprivault} and \cite{penent2022fully} 
to treat the case of systems of fully nonlinear PDEs of the form \eqref{eq:main pde},
 and 
 we apply our algorithm to the incompressible Navier-Stokes equation
\begin{equation}
\nonumber 
\begin{cases}
  \displaystyle
  \partial_t u(t,x) + \nu \Delta u(t,x)
  =
  \nabla p(t,x)
 + ( u \cdot \nabla ) u,
  \medskip
  \\
  \Delta p(t, x)
  = - {\rm div \ \! }{\rm div \ \! } ( u \otimes u ),
  \medskip
  \\
u(T,x) = \phi (x),
\quad (t,x) = (t,x_1, \ldots, x_d) \in [0,T] \times \real^d,
\end{cases}
\end{equation}
with pressure term $p(t,x)=u_0(t,x)$.
This equation is a special case of \eqref{eq:main pde}
obtained by taking $n=d(d+2)$ and $q=d$,
 see Section~\ref{sec:numerical examples},
 and can be rewritten as the divergence-free problem
\begin{equation}
\nonumber 
\begin{cases}
  \displaystyle
  \partial_t u(t,x) + \nu \Delta u(t,x)
  =
  \nabla p(t,x)
 + ( u \cdot \nabla ) u,
  \medskip
  \\
{\rm div \ \! } u = 0,
  \medskip
  \\
u(T,x) = \phi (x),
\quad (t,x) = (t,x_1, \ldots, x_d) \in [0,T] \times \real^d.
\end{cases}
\end{equation}

\medskip

 Probabilistic representations for the solution of the Navier-Stokes equation
 using BSDEs have been considered in e.g.
 \cite{belopolskaya}, \cite{cipriano}, \cite{cruzeiroshamarova}, and 
 Monte Carlo numerical algorithms based on BSDEs have been designed and
implemented in \cite{delbaen2015forward} and \cite{lejay2020forward}.
The BSDE approach is however restricted to first order nonlinear PDE systems
for which we have $\max\limits_{j} \sum\limits_{i=1}^d \alpha^j_i \leq 1$
in \eqref{eq:main pde}, and its numerical implementation
 involves errors from both Monte Carlo estimation
 and time discretization, thus reducing its effectiveness. 
 The Navier-Stokes equation can also be solved in
 the framework of physics-informed neural networks (PINN) using
 the Galerkin method \cite{Rassi19} 
 and high quality solution data usually obtained from an existing solver over
 a given training domain.
 On the other hand, our branching algorithm belongs to the family
     of solvers that do not use existing training data. 
 
\medskip

In Section~\ref{sec:numerical examples} we compare our numerical results to
those obtained in \cite{lejay2020forward} using BSDEs and the Monte Carlo method,
and in \cite{angeli} using finite-difference and finite-element methods.
We note in particular that our method is more stable and much faster
than the BSDE approach of \cite{lejay2020forward}
which has been implemented on a computer cluster with a few tens of cores.

\medskip

We also compare our results to those of \cite{angeli} in which
the 2D Taylor-Green example has been treated by finite-difference and finite-element methods with viscosities $\nu = 10^{-1}$ and $\nu =10^{-2}$, see Section~5 therein.
Although we cannot fully match the speed and precision of state of the art finite element methods, we would like to stress the following points.
     \begin{itemize}
     \item Our neural network approach yields a full functional estimation
       over a whole time-space domain instead of
 pointwise estimates on a grid as in mesh-based methods.
\item
     Monte Carlo estimation provides an intuitive interpretation
     of the solution of partial differential
     equations via the diffusion of heat mechanism,
     as such it makes sense to test their applicability,
     which also opens the door to future applications to the solution 
     of higher dimensional systems. 

        In particular, our branching algorithm overcomes the curse
     of dimensionality because the number of tree branches is not
     sensitive to dimension, see the comments at the end of
     Section~\ref{subsec:feynman kac}. 
       For example, in
       \cite{nguwipenentprivault}, \cite{penent2022fully} 
       this branching method has been applied to PDE examples in dimension $100$,
       which may not be treated using mesh-based methods. 
     \item Our results compare favorably to other
       Monte Carlo algorithms such as \cite{lejay2020forward}, 
       in which computations for a single time step can require several hours. 
     \end{itemize}
 
     In \cite{mmatsumoto}, the {D}eep {G}alerkin {M}ethod (DGM) has been applied to the numerical solution of compressible {N}avier-{S}tokes equations with Reynolds numbers around 1,000, and in \cite{li-yue-zhang-duan},
     the DGM method has been applied  to time-independent Stokes equations. 
     However, we have not encounter applications of the DGM method
     to the incompressible Navier-Stokes equation in the literature,
     including for the Taylor-Green vortex  and the Arnold-Beltrami-Childress flow.
 In Section~\ref{fjkd13} we compare the output of our method
 to that of the deep Galerkin method \cite{sirignano2018dgm},
 see Figures~\ref{fig:dgm1}-\ref{fig:dgm4}.
 We note that the DGM method performs correctly if one reduces the domain of
study from $[0,2\pi]^2$ to $[0,1]^2$ as done in e.g. \cite{li-yue-zhang-duan} 
 for Stokes equations, and uses space-time boundary conditions on $[0,1]^2\times [0,T]$.
On the other hand, we observe that the DGM algorithm loses its accuracy 
when only a condition at terminal time $T$ is used 
as in our method, or when the domain is extended from $[0,1]^2$ to
$[0,2\pi]^2$, see Figures~\ref{fig:dgm2}-\ref{fig:dgm4}.

     \medskip

     Although our method does not use 
     domain boundary conditions, such conditions
     can be taken into account by replacing
     the standard Gaussian kernel by specialized
     heat kernels, see e.g.
     \S~III-4 of \cite{borodin}
     for explicit heat kernel expressions with rectangle
     boundary conditions. 

     \medskip 

 In addition to dealing with the Navier-Stokes equation,
 the framework of Equation~\eqref{eq:main pde} is general
 enough to potentially cover equations of non-Newtonian fluid mechanics
 in which viscosity may depend on the gradient of the solution,
 as, for example, in the non-Newtonian Navier-Stokes equation
$$
 \partial_t u(t,x) +
 \xi \nu | \partial_x u(t, x)|^{\xi - 1} \Delta u(t, x)
  =
  \nabla p(t,x)
 + ( u \cdot \nabla ) u,
 $$
  for a power-law non-Newtonian flow, here in dimension $d=1$,
  $\xi > 0$.

    \medskip

We proceed as follows. In Section~\ref{subsec:feynman kac},
we present the construction of the probability representation
\eqref{eq:feynman kac} for the PDE system \eqref{eq:main pde} with
the corresponding algorithm.
Then, in Section~\ref{sec:main ideas}
we outline the deep branching method
for the estimation of \eqref{eq:feynman kac}.
Then in Section~\ref{sec:numerical examples}
 we apply the deep branching method to the examples of Taylor-Green vortex
 and Arnold-Beltrami-Childress flow, 
 and we present further examples using rotating flows.
 
 \medskip

  The Python codes used for our numerical illustrations
 are available at 
 \url{https://github.com/nguwijy/deep_navier_stokes}.

\subsubsection*{Notation}
We denote by $\nn = \{0, 1, 2, \dots\}$ the set of natural numbers.
We let ${\cal C}^{0, \infty}([0, T] \times \real^d; \real^k)$
be the set of functions $u:[0, T] \times \real^d \to \real^k$
such that $u(t, x)$ is continuous in $t$ for all $x\in \real^d$,
and infinitely $x$-differentiable for all $t\in [0,T]$.
For a vector $x = (x_1, \ldots, x_d)^\top \in \real^d$,
we let $\abs{x} = \sum\limits_{i=1}^d \abs{x_i}$
and 
$\bm{1}_p$ be the indicator vector made of $1$ at position $p\in \{1,\ldots , d\}$, and $0$ elsewhere.
As in \cite{constantine}, for use in the multivariate Fa\`a di Bruno formula
we also define a linear order $\prec$ on $\real^d$ such that
$(k_1, \dots, k_d) = k  \prec l = (l_1, \dots, l_d)$
if one of the following holds:
\begin{enumerate}[label=\emph{\roman*})]
    \item $\abs{k} < \abs{l}$;
    \item $\abs{k} = \abs{l}$ and $k_1 < l_1$;
    \item $\abs{k} = \abs{l}$, $k_1 = l_1, \dots k_i = l_i$,
            and $k_{i+1} < l_{i+1}$ for some $1 \leq i < d$.
\end{enumerate}
\noindent
Given $\mu \in \nn^d$, $f \in {\cal C}^\infty ( \real^n)$
and $v \in {\cal C}^{0, \infty}([0, T] \times \real^d;\real^n)$,
 we will use the multivariate Fa\`a di Bruno formula
 \begin{equation}
 \label{eq:fdb}
 \partial_\mu ( f (v(t, x)) ) =
 \left(\prod_{i = 1}^d \mu_i! \right)
    \sum\limits_{\substack{1 \leq \lambda_1 + \cdots + \lambda_n \leq \abs{\mu} \\
                           1 \leq s \leq \abs{\mu}}}
   \partial_{\lambda} f (t,x)
    \hskip-0.64cm
    \sum\limits_{\substack{1 \leq \abs{k_1}, \dots, \abs{k_s}, \
                           0 \prec l^1 \prec \cdots \prec l^s \\
                           k^i_1 + \cdots + k^i_s = \lambda_i, \
                           i = 1, \dots, n \\
                           \abs{k_1}l_j^1 + \cdots + \abs{k_s}l_j^s = \mu_j, \
                           j = 1, \dots, d
                           }}
    \prod_{\substack{1 \leq i \leq n \\
                     1 \leq r \leq s
                    }}
\frac{(\partial_{l^r } v_i (t,x))^{k_r^i}}{k_r^i!
                \left(l_1^r! \cdots l_d^r!\right)^{k_r^i}},
 \end{equation}
 $x=(x_1,\ldots , x_d) \in \real^d$,
 see Theorem~2.1 in \cite{constantine}.
\section{Fully nonlinear Feynman-Kac formula}
\label{subsec:feynman kac}
In this section we extend the construction of
\cite{penent2022fully}, \cite{nguwipenentprivault} to the case of
systems of fully nonlinear coupled parabolic
 and Poisson equations PDEs of the form \eqref{eq:main pde}.
 For this, we rewrite \eqref{eq:main pde} in integral form
 for $i = 1,\ldots , d$ as
\begin{numcases}{}
  \label{eq:integral u0-1} 
  \displaystyle
    u_0(t,x) =
    \frac{\Gamma (d/2)}{2 \pi^{d/2}}
    \int_{\real^d}
  \frac{N(y)}{\abs{y}^d}
 f_0\big(
            \partial_{\bar{\alpha}^{q+1}}u_{\beta_{q+1}}(t,x+y)
        ,
        \ldots ,
                   \partial_{\bar{\alpha}^n}u_{\beta_n}(t,x+y)
        \big)
        dy, 
  \bigskip
 \\
   \label{eq:integral u0}
  \displaystyle
 u_i (t,x) =  \int_{\real^d} \varphi_{2\nu}( T-t,y-x)\phi_i(y)dy 
   \\
   \nonumber
   \displaystyle
 + \int_t^T \int_{\real^d} \varphi_{2\nu} ( s-t ,y-x)
 f_i\big(
    \partial_{\bar{\alpha}^1}u_0 (s,y)
    ,
    \ldots ,
    \partial_{\bar{\alpha}^q }u_0 (s,y)
    ,
        \partial_{\bar{\alpha}^{q+1}}u_{\beta_{q+1}}(s,y)
        ,
        \ldots ,
                   \partial_{\bar{\alpha}^n}u_{\beta_n}(s,y)
        \big)
      dy ds,
  \bigskip
 \\
  \nonumber
  u_i(T,x) = \phi_i (x),
\quad (t,x) = (t,x_1, \ldots, x_d) \in [0,T] \times \real^d,
\ \
i = 0,1,\ldots , d, \
 (t,x) \in [0,T]\times \real^d,
\medskip
\end{numcases}
\addtocounter{equation}{1}
under appropriate integrability condition
as in e.g. Lemma~1.6 in \cite{majda2002vorticity}, where
 $\varphi_{\sigma^2} (t,x) = \re^{- x^2 / ( 2\sigma^2 t) } / {\sqrt{2\pi \sigma^2 t} }$,
$$
\phi_0(x) := u_0(T, x) =
\frac{\Gamma (d/2)}{2 \pi^{d/2}}
\int_{\real^d}
 \frac{N(y)}{  \abs{y}^d}
 f_0\big(
          \partial_{\bar{\alpha}^{q+1}}\phi_{\beta_{q+1}}(x+y)
        ,
        \ldots ,
                   \partial_{\bar{\alpha}^n}\phi_{\beta_n}(x+y)
        \big)
    dy,
$$
 $\Gamma (y) := \int_0^\infty x^{z-1} \re^{-x} dx$ is the Gamma function,
 and $N(y)$ is the Poisson kernel
$$
N(y) =
\begin{cases}
  \displaystyle
  \abs{y}^2 \log\abs{y},     & d = 2,
    \medskip
    \\
  \displaystyle
    \frac{\abs{y}^2}{2-d},      & d \ge 3, \quad y\in \real^d.
\end{cases}
$$
 Our fully nonlinear Feynman-Kac formula 
 relies on the
 construction of a branching coding tree, based on the definition of
 a set $\mathcal{C}$ of codes and its associated branching mechanism $\mathcal{M}$.
In what follows, for any function $g:\real^n \to \real$,
we let $g^*$ be the operator mapping
${\cal C}^{0, \infty}([0, T] \times \real^d)$
to
${\cal C}^{0, \infty}([0, T] \times \real^d)$
and defined by
\begin{eqnarray*}
g^*(u)(t, x) & := &
    g \big(
    \partial_{\bar{\alpha}^1}u_0 (t,x)
    ,
    \ldots ,
    \partial_{\bar{\alpha}^q }u_0 (t,x)
    ,
        \partial_{\bar{\alpha}^{q+1}}u_{\beta_{q+1}}(t,x)
        ,
        \ldots ,
                   \partial_{\bar{\alpha}^n}u_{\beta_n}(t,x)
                   \big)
                   \\
                    & = &
                   g\big(\partial_{\bar{\alpha}^1}u_{\beta_1}(t,x),\ldots , \partial_{\bar{\alpha}^n}u_{\beta_n}(t,x)\big),
\end{eqnarray*}
with $\beta_1=\cdots =\beta_q=0$,
from which \eqref{eq:main pde} can be rewritten as
 $$
  \partial_t u_i(t,x) + \nu \Delta u_i(t,x)
  +
  f_i^*(u)(t, x) = 0,
  \qquad
  i = 1,\ldots , d,
  $$
 $(t,x) = (t,x_1, \ldots, x_d) \in [0,T] \times \real^d$.
  In the sequel, for $\lambda = \left(\lambda_1, \dots, \lambda_n\right) \in \nn^n$
we let
$$
\partial_{\lambda} f_i(x_1,\ldots , x_n)
=
\frac{\partial^{\lambda_1}}{\partial x_1^{\lambda_1}}
\cdots
\frac{\partial^{\lambda_n}}{\partial x_n^{\lambda_n}}
f_i(x_1,\ldots , x_n), \qquad (x_1,\ldots , x_n)\in \real^n,
$$
$$
\partial_{\lambda} f_0(x_{q+1},\ldots , x_n)
=
\frac{\partial^{\lambda_1}}{\partial x_1^{\lambda_1}}
\cdots
\frac{\partial^{\lambda_n}}{\partial x_n^{\lambda_n}}
f_0(x_{q+1},\ldots , x_n), \qquad (x_1,\ldots , x_n)\in \real^n.
$$

\vspace{-0.3cm}

\begin{definition}
 We let $\mathcal{C}$ denote the set of
 operators from
 ${\cal C}^{0,\infty} ([0,T]\times \real^d; \real^{d+1} )$
 to \\
 ${\cal C}^{0,\infty} ([0,T]\times \real^d; \real )$
 called \textit{codes}, and defined as
\begin{equation} 
  \label{fjkl214}
  \mathcal{C} := \left\{
           {\rm Id}_i, \
           (a \partial_{\lambda} f)^*, \
           (a \partial_{\mu}, i), \
           (\partial_{\mu}, -1) \ : \
           \lambda \in \mathbb{N}^n, \
           \mu \in \mathbb{N}^d, \
           a \in \real, \
                      i = 0, \ldots , d
       \right\}.
\end{equation} 
\end{definition}
\noindent
 The codes $c$ in $\mathcal{C}$ are operators acting on
$(u_0, u_1, \ldots, u_d) = u
 \in {\cal C}^{0,\infty} \left([0,T]\times \real^d; \real^{d+1}\right)$
 as
\begin{equation*}
    c(u)(t, x) =
    \begin{cases}
      \displaystyle
      u_i(t, x),
                & \text{if } c = {\rm Id}_i,
      \medskip
      \\
      a \partial_\lambda f\left(
        \partial_{\bar{\alpha}^1}u_{\beta_1}(t, x)
, \ldots ,
        \partial_{\bar{\alpha}^n}u_{\beta_n}(t, x)
      \right),
                & \text{if } c = (a \partial_\lambda f)^*,
      \medskip
      \\
      a \partial_\mu u_i(t, x),
                & \text{if } c = (a\partial_\mu, i),
      \medskip
      \\
      \partial_\mu (\partial_t + \nu \Delta) u_0(t, x),
                & \text{if } c = (\partial_\mu, -1).
    \end{cases}
\end{equation*}
 In the branching algorithm implementation, 
 the quantities $\partial_\lambda f\left(
        \partial_{\bar{\alpha}^1}u_{\beta_1}(t, x) , \ldots , \partial_{\bar{\alpha}^n}u_{\beta_n}(t, x)
        \right)$
        and $\partial_\mu u_i(t, x)$
        will be computed by estimating $u_i(t,x)$ recursively
        using the integral form of \eqref{eq:main pde} for $i\in \{1,2,\ldots , d\}$.
        No such recursion is needed for $\partial_\mu u_0(t, x)$
        and $\partial_\mu (\partial_t + \nu \Delta) u_0(t, x)$ which will be
        computed by solving the corresponding Poisson equation using
        integral expressions, see
        \eqref{eq:poisson integration calculation}
        and
        \eqref{eq:poisson integration calculation-2}
        in appendix.

        \medskip

        This recursion will be implemented
        using the {\em branching mechanism} $\mathcal{M}$ defined below,
        which is based on the multivariate Fa\`a di Bruno formula \eqref{eq:fdb}.
 For the description and implementation of the algorithm
        we will enumerate the terms appearing
        in \eqref{eq:fdb}
        applied to the index set $\mu \in \mathbb{N}^n$ and function
        $f$ on $\real^n$
 using the set
${\rm fdb}(\mu, f, (c_1, \dots, c_m))$ of code sequences
  defined as
\begin{align*}
& {\rm fdb}(\mu, f, (c_1, \dots, c_m))
\\
& :=
\bigcup_{
     \footnotesize \substack{1 \leq s \leq \abs{\mu}, \
               1 \leq \lambda_1 + \cdots + \lambda_n \leq \abs{\mu} \\
               1 \leq \abs{k_1}, \dots, \abs{k_s}, \
               0 \prec l^1 \prec \cdots \prec l^s \\
               k^i_1 + \cdots + k^i_s = \lambda_i, \
               i = 1, \dots, n \\
               \abs{k_1}l_j^1 + \cdots + \abs{k_s}l_j^s = \mu_j, \
               j = 1, \dots, d
                }}
\left\{
(c_1, \dots, c_m) \ \bigcup
\left(
\frac{\prod\limits_{j = 1}^d \mu_j! }
     {\prod\limits_{\footnotesize \substack{1 \leq i \leq n \\
                     1 \leq r \leq s}}
      k_r^i! \left(l_1^r! \cdots l_d^r!\right)^{k_r^i}
     }
     (\partial_{\lambda} f)^*
     \right)
     \right.
\\
 & \qquad \qquad \qquad \qquad \qquad \qquad \qquad \qquad \qquad \qquad \left.
     \bigcup
     \bigcup_{\footnotesize \substack{1 \leq i \leq n \\
                     1 \leq r \leq s
                    }}
     \big(
     \underbrace{
         (\partial_{l^r + \bar{\alpha}^i}, \beta_i)
       ,
       \ldots
       ,
       (\partial_{l^r + \bar{\alpha}^i}, \beta_i)
       }_{k_r^i~{\rm times}}
     \big)
\right\},
\end{align*}
where $(c_1, \dots, c_m)$ is any sequence of codes in
$\mathcal{C}$ and we use the notation
$$(a_1 , \ldots , a_n) \cup (b_1 , \ldots , b_m)
:=
(a_1 , \ldots , a_n, b_1 , \ldots , b_m)
$$
and
$(a_1 , \ldots , a_n) \cup \emptyset = (a_1 , \ldots , a_n)$
for any sequences
$(a_1 , \ldots , a_n)$, $(b_1 , \ldots , b_m)$.  
The next definition provides a way to enumerate the terms appearing
in \eqref{eq:fdb} and in \eqref{eq:poisson integration calculation}-\eqref{eq:poisson integration calculation-2} below.
\begin{definition}
\label{fjklds114} 
 We define a {\em mechanism} $\mathcal{M}$
 that maps any code $c$ in $\mathcal{C}$ to a set ${\cal M}(c)$ of code
 sequences, by letting
\begin{align*}
  & \mathcal{M} ( {\rm Id}_i ) := \{ f_i^* \},
  \quad            i = 0,1,\ldots , d,
  \medskip
  \\
 & {\cal M}\left((\partial_{\mu}, i)\right)
:= {\rm fdb}\left(\mu, f_i, \emptyset \right),
  \qquad
  \mu \in \mathbb{N}^n,
  \quad            i = 0,1,\ldots , d,
  \medskip
\\
  &
  {\cal M}(g^*)
 :=
 \bigcup_{\footnotesize \substack{q < r \leq n}}
      {\rm fdb}\left(\bar{\alpha}^r, f_{\beta_r}, ((\partial_{\bm{1}_r} g)^*) \right)
 \ \bigcup
\bigcup_{
     \footnotesize \substack{i, j = 1, \dots, n \\
               k = 1, \dots, d
                }}
\left\{
\left(-\nu (\partial_{\bm{1}_i + \bm{1}_j} g)^*,
    (\partial_{\bar{\alpha}^i + \bm{1}_k}, \beta_i),
(\partial_{\bar{\alpha}^j + \bm{1}_k}, \beta_j)
\right)
\right\}
\medskip
\\
& \qquad \qquad \qquad
\bigcup
\bigcup_{1 \leq r \leq q}
\left\{
\bigl(
    -(\partial_{\bm{1}_r} g)^*,
    (\partial_{\bar{\alpha}^r}, -1)
\bigr)\right\},
\qquad g^* \in \mathcal{C},
\end{align*}
 and
\begin{align*}
& {\cal M}((\partial_\mu, -1))
  \\
   &
  :=
  \bigcup_{
    {
      {
        {
          q < i, j \leq n
          \atop k = 1, \dots, d}
    \atop 0 \leq \ell_p \leq \gamma_p \leq \mu_p }
      \atop p = 1, \dots, d
    }
    }
{\rm fdb}\left(
    \gamma,
    \partial_{\bm{1}_i + \bm{1}_j} f_0,
    \left(
    \left(
        \left(
        \nu
        \prod\limits_{r = 1}^d
        {\mu_r \choose \ell_r}
        {\ell_r \choose \gamma_r}
        \right)
    \partial_{\mu - \ell + \bar{\alpha}^i + \bm{1}_k}, \beta_i
    \right),
    \left(
            \partial_{\ell - \gamma + \bar{\alpha}^j + \bm{1}_k}, \beta_j
    \right)
    \right)
\right)
\\
&
\bigcup \
\bigcup_{
  {
    {
      q < i \leq n
       \atop
      {0 \leq \ell_p \leq \mu_p
        \atop p = 1, \dots, d
        }
    }
  }
  }
\bigcup_{
     \footnotesize \substack{1 \leq s \leq \abs{\ell}, \
               1 \leq \lambda_1 + \cdots + \lambda_n \leq \abs{\ell} \\
               1 \leq \abs{k_1}, \dots, \abs{k_s}, \
               0 \prec l^1 \prec \cdots \prec l^s \\
               k^i_1 + \cdots + k^i_s = \lambda_i, \
               i = 1, \dots, n \\
               \abs{k_1}l_j^1 + \cdots + \abs{k_s}l_j^s = \ell_j, \
               j = 1, \dots, d
                }}
\\
& 
{\rm fdb}\left(
    \mu - \ell + \bar{\alpha}^i, f_{\beta_i},
    \left(
    - \frac{
        \prod\limits_{j = 1}^d
        \frac{\mu_j!}{(\mu_j-\ell_j)!}
            }
         {\prod\limits_{\footnotesize \substack{1 \leq i \leq n \\
                         1 \leq r \leq s}}
          k_r^i! \left(l_1^r! \cdots l_d^r!\right)^{k_r^i}
         }
         (\partial_{\lambda + \bm{1}_i} f_0)^*
         \right)
     \bigcup
     \bigcup_{\footnotesize \substack{1 \leq i \leq n \\
                     1 \leq r \leq s
                    }}
     \big(
     \underbrace{
         (\partial_{l^r + \bar{\alpha}^i}, \beta_i)
       ,
       \ldots
       ,
       (\partial_{l^r + \bar{\alpha}^i}, \beta_i)
       }_{k_r^i~{\rm times}}
     \big)
    \right).
\end{align*}
\end{definition}
In order to motivate the construction of the mechanism ${\cal M}$,
we note that
\begin{itemize}
\item ${\cal M}\left((\partial_{\mu}, i)\right)$, 
  $i=1,\ldots , d$,
  is used to model the Fa\`a di Bruno formula
  \eqref{eq:fdb} 
 via ${\rm fdb}\left(\mu, f_i, \emptyset \right)$, 
\item $\mathcal{M} ( {\rm Id}_0 )$
  is used to model the 
  Poisson integral equation \eqref{eq:integral u0-1},
\item $\mathcal{M} ( {\rm Id}_i )$,
         $i=1,\ldots ,d$, 
  is used to model the integral equation \eqref{eq:integral u0},
\item 
  ${\cal M}(g^*)$ is used to model the integral
  equation \eqref{fjkld133}, 
\item
${\cal M}((\partial_\mu, 0))$ 
 is used to model the Poisson integral equation \eqref{eq:poisson integration calculation}, 
\item ${\cal M}((\partial_\mu, -1))$
  is used to model the integral equation \eqref{jkldf09}. 
\end{itemize}
\subsubsection*{Example - semilinear PDEs}
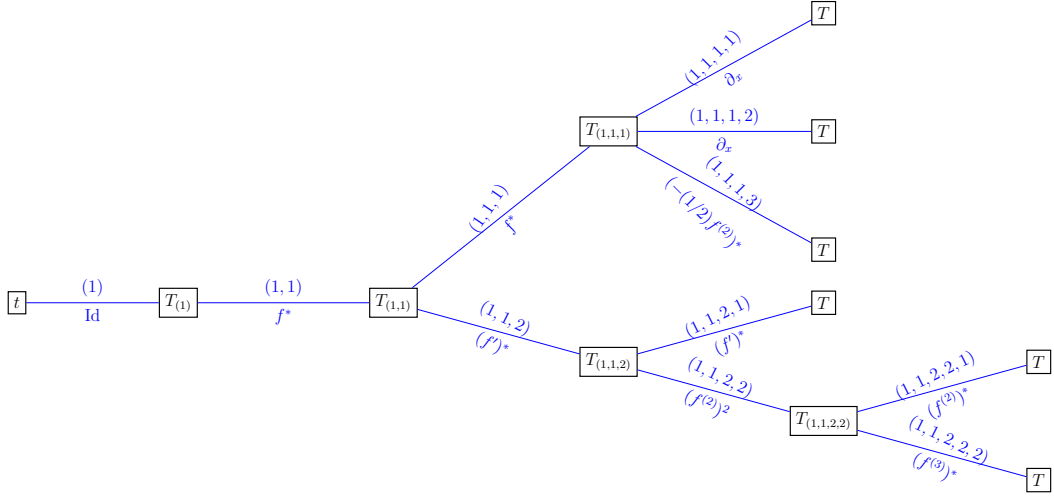
\begin{figure} 
\centering
\resizebox{0.85\textwidth}{!}{
\begin{tikzpicture}[yscale = 0.55,scale=1.3,grow=right, sloped,color=blue]
\node[rectangle,draw,black,text=black,thick]{$t$}
            child {
              node[rectangle,draw,black,text=black,thick] {$T_{(1)}$}
              child {
                node[rectangle,draw,black,text=black,thick] {$T_{(1,1)}$}
                child{
                node[rectangle,draw,black,text=black,thick]{$T_{(1,1,2)}$}
                    child{
                node[rectangle,draw,black,text=black,thick]{$T_{(1,1,2,2)}$}
                child{
                    node[rectangle,draw,black,text=black,thick]{$T$}
                    edge from parent
                    node[above]{$(1,1,2,2,2)$}
                    node[below]{$(f^{(3)})^*$}
                    }
                    child{
                    node[rectangle,draw,black,text=black,thick]{$T$}
                    edge from parent
                    node[above]{$(1,1,2,2,1)$}
                    node[below]{$( f^{(2)})^*$}
                    }
                    edge from parent
                node[above]{$(1,1,2,2)$}
                node[below]{$( f^{(2)})^2$}
                }
                    child{
                    node[rectangle,draw,black,text=black,thick]{$T$}
                    edge from parent
                    node[above]{$(1,1,2,1)$}
                    node[below]{$(f')^*$}
                    }
                edge from parent
                node[above]{$(1,1,2)$}
                node[below]{$(f')^*$}
                }
                child{
                node[above=10pt,rectangle,draw,black,text=black,thick,yshift=2cm]{$T_{(1,1,1)}$}
                    child{
                    node[rectangle,draw,black,text=black,thick]{$T$}
                    edge from parent
                    node[above]{$(1,1,1,3)$}
                    node[below,yshift=-0.3cm]{$(-(1/2) f^{(2)})^*$}
                    }
                    child{
                    node[rectangle,draw,black,text=black,thick]{$T$}
                    edge from parent
                    node[above]{$(1,1,1,2)$}
                    node[below]{$\partial_x$}
                    }
                    child{
                    node[rectangle,draw,black,text=black,thick]{$T$}
                    edge from parent
                    node[above]{$(1,1,1,1)$}
                    node[below]{$\partial_x$}
                    }
edge from parent
                node[above]{$(1,1,1)$}
                node[below]{$f^*$}
                                }
                edge from parent
                node[above] {$(1,1)$}
                node[below]  {$f^*$}
            }
                edge from parent
                node[above] {$(1)$}
                node[below]  {${\rm Id}$}
            };
\end{tikzpicture}
}
\vspace{-0.2cm}
\caption{Sample coding tree.}
\label{fjkldsf-f}
\end{figure}

 To illustrate our method,
 consider the simpler case of a semilinear PDE of the form
\begin{equation}
   \label{eq:112}
\begin{cases}
  \displaystyle
  \partial_t u(t,x) + \frac{1}{2}\partial_x^2 u(t,x) + f(u(t,x)) = 0
  \medskip
  \\
u(T,x) = \phi(x), \qquad (t,x) \in [0,T] \times \real,
\end{cases}
\end{equation}
 in dimension $d=1$, with the integral formulation
\begin{equation}
\label{fjkldsf}
 u(t,x) = \int_{-\infty}^\infty \varphi (T-t,y-x) \phi(y) dy
 + \int_t^T \int_{-\infty}^\infty \varphi (s-t,y-x) f(u(s,y)) dy ds.
\end{equation}
 In order to iterate \eqref{fjkldsf} into a tree-based recursion,
 we will find a PDE satisfied by $f(u(s,y))$ by differentiating
\begin{equation*}
\begin{split}
  \partial_s f(u(s,y)) + \frac{1}{2}\partial_y^2 f(u(s,y))
  &=
  \left( \partial_s u (s,y)+ \frac{1}{2} \partial_y^2 u(s,y) \right)
  f'(u(s,y))
  +  \frac{1}{2} (\partial_y u(s,y))^2 f''(u(s,y))
  \\
&= -f(u(s,y)) f'(u(s,y)) + \frac{1}{2}(\partial_y u(s,y))^2 f''(u(s,y)),
\end{split}
\end{equation*}
 showing that $f(u(s,y))$ satisfies  the integral equation
 \begin{eqnarray}
   \label{jfkdls32}
   \lefteqn{ 
  f(u(s,y)) = \int_{-\infty}^\infty \varphi (T-s,z-x) f( \phi(z) ) dz
   }
   \\
   \nonumber
   & & 
   + \int_s^T \int_{-\infty}^\infty \varphi ( w -s ,z-x)
 \left( f(u( w ,z)) f'(u( w ,z)) - \frac{1}{2}(\partial_z u( w ,z))^2 f''(u( w ,z)) \right) dz d w,
\end{eqnarray}
 More generally, we use \eqref{fjkldsf} and \eqref{jfkdls32} 
 to expand $u(t,x)$ and $a f^{(k)} ( u(t,x) )$
 into a consistent set of equations
 which are suitable for a recursive estimation of $u(t,x)$, as
$$
  \left\{
  \begin{array}{l}
    \displaystyle
    u(t,x) = \int_{-\infty}^\infty \varphi (T-t,y-x) \phi(y) dy + \int_t^T \int_{-\infty}^\infty \varphi (s-t,y-x) f(u(s,y)) dy ds
    \medskip
    \\
    \displaystyle
  a f^{(k)} ( u(t,x) )= \int_{-\infty}^\infty \varphi (T-t,y-x) a f^{(k)} ( \phi (y) ) dy
\medskip
    \\
 \quad
    \displaystyle
 + \int_t^T \int_{-\infty}^\infty \varphi (s-t,y-x)
\left( af ( u(s,y)) f^{(k+1)} ( u(s,y)) - \frac{a}{2}(\partial_y u (s,y))^2 f^{(k+2)} ( u (s,y)) \right) dy ds
  \end{array}
  \right.
  $$
  $a \in \real\setminus \{0\}$, $k \geq 0$, and we expand $\partial_x u(t,x)$ as 
$$
  \partial_x u(t,x) =  \int_{-\infty}^\infty \varphi (T-t,y-x) \partial_x \phi(y) dy + \int_t^T \int_{-\infty}^\infty \varphi (s-t,y-x) f' (  u (s,y) ) \partial_y u (s,y) dy  ds.
  $$
   \noindent
 In this case, the set of $\eufrak{C}$ codes in \eqref{fjkl214} reads 
$$
 \eufrak{C} := \big\{
           {\rm Id}, \ \partial_x, \ a f^{(k)}, \ a \in \real\setminus \{0\}, \ k \in \mathbb{N} \big\},
$$
           and the branching mechanism $\mathcal{M}$ in Definition~\ref{fjklds114}
           is given by
                 \begin{equation}
                   \label{jfklds}
 {\cal M}({\rm Id}) := \{f^* \},
 \ \
     {\cal M}(g^*) := \left\{\big(f^*,(\partial_{z_0}g)^*\big)
     ;
      \left(\partial_x,\partial_x, -\frac{1}{2}(\partial_{z_0}^2 g)^* \right)\right\},
 \ \
 {\cal M}(  \partial_x ) := \big\{ \big( (\partial_{z_0}f)^* ,\partial_x \big) \big\},
\end{equation}
                 for $g\in {\cal C}^\infty ( \real^{n+1})$
                  of the form $g = a \partial_{z_0}^k f$, $a \in \real\setminus \{0\}$, $k \geq 0$.
  Figure~\ref{fjkldsf-f} presents a sample of the random coded tree
  $\mathcal{T}_{t,x,{\rm Id}}$ started from $c = {\rm Id}$
  for a semilinear PDE of the form \eqref{eq:112}.

  \medskip

  \tikzstyle{level 1}=[level distance=3cm, sibling distance=1cm]
\tikzstyle{level 2}=[level distance=4cm, sibling distance=4cm]

\subsubsection*{Implementation} 
 The probabilistic representation of PDE solutions
 will be implemented using the functional $\mathcal{H}({t, x, c})$
 constructed in Algorithm~\ref{alg:coding tree} below
 along a random coding tree started at
$(t,x,c)\in [0,T]\times \real^d \times {\cal C}$.
 We consider two probability density functions (PDF)
 $\rho, \widetilde{\rho}: \real_+ \to (0, \infty)$ on $\real_+$, and denote by
 $\widebar{F}$ the tail distribution function of $\rho$,
 and let $\mathcal{N}(0, \sigma^2 {\rm I}_d )$ denote the $d$-dimensional
 centered normal distribution with variance $\sigma^2$ and independent
 components.

 \begin{algorithm}[H]
  \caption{Coding tree algorithm TREE$(t, x, c)$.}
  \label{alg:coding tree}
    \begin{algorithmic}
        \Require $t \in [0, T]$, $x \in \real^d$, $c \in \mathcal{C}$
        \Ensure $\mathcal{H}({t, x, c}) \in \real$
        \State $\mathcal{H}({t, x, c}) \gets 1$
        \State $\tau \gets$ a random variable drawn from the distribution of $\rho$
        \State $\widetilde{\tau} \gets$ a random variable drawn from the distribution of $\widetilde{\rho}$
        \If{$t + \tau > T$}
            \State $W_{2 \nu (T - t)} \gets$ a random vector drawn from $\mathcal{N}(0, 2 \nu (T - t))$
            \State $\mathcal{H}({t, x, c}) \gets \mathcal{H}({t, x, c}) \times c(u)(T, x + W_{2 \nu (T - t)})
            /  \widebar{F}(T-t)$
        \ElsIf{$c \in \left\{ (a\partial_\mu, i): a \in \real, \mu \in \nn^d, i = 0, -1 \right\}$}
            \State $W_{\widetilde{\tau}} \gets$ a random vector drawn from $\mathcal{N}(0, \widetilde{\tau})$
            \State $r_c \gets$ the size of the mechanism set $\mathcal{M}(c)$
            \State $I_c \gets$ a random element drawn uniformly from $\mathcal{M}(c)$
            \State $\mathcal{H}({t, x, c}) \gets \mathcal{H}({t, x, c}) \times N(W_{\widetilde{\tau}}) \times r_c \times (2 \widetilde{\tau} \ \! \widetilde{\rho}(\widetilde{\tau}))^{-1}$
            \ForAll{$cc \in I_c$}
                \State $\mathcal{H}({t, x, c}) \gets \mathcal{H}({t, x, c}) \times \text{TREE}(t, x + W_{\widetilde{\tau}}, cc)$
            \EndFor
        \Else
            \State $W_{2\nu \tau} \gets$ a random vector drawn from $\mathcal{N}(0, 2 \nu \tau)$
            \State $r_c \gets$ the size of the mechanism set $\mathcal{M}(c)$
            \State $I_c \gets$ a random element drawn uniformly from $\mathcal{M}(c)$
            \State $\mathcal{H}({t, x, c}) \gets \mathcal{H}({t, x, c}) \times r_c \times \rho^{-1}(\tau)$
            \ForAll{$cc \in I_c$}
                \State $\mathcal{H}({t, x, c}) \gets \mathcal{H}({t, x, c}) \times \text{TREE}(t + \tau, x + W_{2\nu \tau}, cc)$
            \EndFor
        \EndIf
    \end{algorithmic}
\end{algorithm}
\noindent
 As in Theorem~3.2 in \cite{penent2022fully},
 the following Feynman-Kac type identity
 holds for the solution of \eqref{eq:main pde} holds
under suitable integrability conditions
on ${\cal H} ({t, x, {\rm Id}_i} )$
and smoothness assumptions on the coefficients of
\eqref{eq:main pde}.
\begin{prop}
\label{p1}
Let $T>0$ such that $\E [ | {\cal H} (t, x, c ) |]< \infty$,
 $c\in \mathcal{C}$, $(t,x)\in [0,T]\times \real^d$,
  and consider the system of equations
\begin{equation}
 \label{s1--}
    \begin{cases}
    \displaystyle
c(u)(t,x) =
  \sum\limits_{Z \in \mathcal{M}(c)}
  \int_0^\infty
  \int_{\real^d}
  \varphi_1 (s, y)
  \frac{N(y)}{2s}
  \prod\limits_{z \in Z}
  z(u)(t, x + y)
  dy
  ds,
  \quad
  c = ( \partial_{\mu}, 0), \ c = (\partial_{\mu}, -1),
  \\
 (t,x)\in [0,T]\times \real, \mbox{ and}
\\
  \displaystyle
 c(u)(t,x) = \int_{\real^d} \varphi_{2\nu} (T-t,y-x) c ( u)(T,y) dy
+
\sum_{Z \in \mathcal{M}(c)}
\int_t^T \int_{\real^d} \varphi_{2\nu} (s-t,y-x)
\prod_{z \in Z}  z(u)(s,y) dy ds,
\\
\mbox{for all remaining codes } c\in \mathcal{C},
 \ (t,x)\in [0,T]\times \real.
    \end{cases}
    \end{equation}
If the solution of the above system is unique,
 then the solution of \eqref{eq:main pde} admits the probabilistic representation
\begin{equation}
    \label{eq:feynman kac}
    u_i(t, x) = \E \left[ {\cal H} \left({t, x, {\rm Id}_i} \right) \right],
 \qquad (t,x)\in [0,T]\times \real^d,
  \end{equation}
 $i = 0,1,\ldots , d$.
\end{prop}
The proof of Proposition~\ref{p1} is given in appendix.
It proceeds as in the proof of Theorem~3.2 of \cite{penent2022fully},
by showing that
$$
c(u)(t,x) = \E [ \mathcal{H}(t,x,c) ], \qquad
(t,x) \in [0,T] \times \real,
$$
 for all codes $c \in \mathcal{C}$, which implies
 \eqref{eq:feynman kac} by taking $c = {\rm Id}_i$, $i=1,\ldots , d$.

 \medskip

 In numerical applications, 
the expected value $\E [ \mathcal{H}(t,x,c) ]$
in Proposition~\ref{p1} is estimated as the average
 $$
 \frac{1}{N}\sum_{k=1}^N \mathcal{H}(t,x,c)^{(k)}
 $$
 where
 $\mathcal{H}(t,x,c)^{(1)}, \ldots , \mathcal{H}(t,x,c)^{(N)}$
 are independent samples of $\mathcal{H}(t,x,c)$.
 In this case, the error on the estimate of
 $\E [ \mathcal{H}(t,x,c) ]$
 from the Monte Carlo method can be estimated as the standard deviation
$$
\left( \E \left[ \left( \E [ \mathcal{H}(t,x,c) ]
   - \frac{1}{N}\sum_{k=1}^N \mathcal{H}(t,x,c)^{(k)}
   \right)^2 \right]
\right)^{1/2} =
 \frac{1}{\sqrt{N}} \sqrt{\Var [ \mathcal{H}(t,x,c) ]}.
$$
      The main tunable parameter in the stochastic branching algorithm
      is the distribution $\rho$ of the random branching time $\tau$.
      Higher mean branching times result into shorter trees on average,
      therefore requiring a higher number of Monte Carlo samples
      in order to achieve the same precision level. 
      For example, in the case of an exponentially distributed branching time
      with parameter $\lambda$, the average depth of binary branching trees
      until time $t>0$ is of order
      $e^{\lambda t}$, see e.g. \S~4 of \cite{penent4}. 

      \medskip

      Overall, the impact of dimension $d$
      is on the number of sequences in the mechanism ${\cal M}(c)$, i.e. 
    on the number of possible ways of branching. 
    On the other hand, the complexity of the algorithm is
    determined by the number of branches at each branching time, 
    i.e. on the lengths of coding sequences, which do not 
    depend on the dimension $d$. 
    As a result, the complexity of our method has polynomial growth as a (small) power of the dimension $d$, mostly due to the use of $d+1$ coding trees in the algorithm.

\section{Deep branching solver}
\label{sec:main ideas}
Instead of evaluating \eqref{eq:feynman kac} at a given point
$(t,x) \in [0,T]\times \real^d$, we use the $L^2$-minimality property of expectation to
obtain a functional estimation of
$u = (u_1, \ldots , u_d )$
as $u(\cdot , \cdot ) = v^*(\cdot , \cdot )$
 on the support of a random vector $(\zeta , X)$ on $[0,T] \times \real^d$
 such that $\mathcal{H} ({\zeta , X, {\rm Id}_i }) \in L^2$, where
\begin{equation}
\label{eq:L2 minimality}
v^* =
\argmin\limits_{\left\{ v:\real^{d+1} \to \real^d \ \! : \ \! v(\zeta , X) \in L^2\right\}}
\sum\limits_{i = 1}^d
\E \left[
    \left(\mathcal{H}\left({\zeta , X, {\rm Id}_i }\right)
            - v_i(\zeta , X)\right)^2 \right].
\end{equation}
To evaluate \eqref{eq:feynman kac} on $[0, T] \times \Omega$,
where $\Omega$ is a bounded domain of $\real^d$,
we can choose $(\zeta , X)$ to be a uniform random vector on $[0, T] \times \Omega$.

\medskip

In order to implement the deep learning approximation,
we parametrize $v(\cdot , \cdot )$
using a functional space described below.
Given $\sigma:\real \to \real$ an activation function
such as $\sigma_{\rm ReLU}(x) := \max (0,x)$,
$\sigma_{\tanh}(x) := \tanh(x)$,
$\sigma_{\rm Id}(x) := x$,
 we define the set of layer functions $\mathbb{L}^\sigma_{d_1,d_2}$ by
\begin{equation}
    \label{eq:one layer}
    \mathbb{L}^\sigma_{d_1,d_2} :=
    \bigl\{
        L:\real^{d_1} \to \real^{d_2} \ : \ L(x) = \sigma(Wx + b),
  \ x \in \real^{d_1}, \ W \in \real^{d_2 \times d_1}, \ b \in \real^{d_2}
    \bigr\},
\end{equation}
where $d_1 \geq 1$ is the input dimension,
$d_2 \geq 1$ is the output dimension,
and the activation function
$\sigma$ is applied component-wise to $Wx + b$.
Similarly, when the input dimension and the output dimension are the same,
we define the set of residual layer functions
$\mathbb{L}^{\rho, \rm res}_d$ by
\begin{equation}
    \label{eq:one layer resnet}
    \mathbb{L}^{\sigma, \rm res}_d :=
    \bigl\{
        L:\real^d \to \real^d \ : \ L(x) = x + \sigma(Wx + b),
  \ x \in \real^d, \ W \in \real^{d \times d}, \ b \in \real^d
    \bigr\},
\end{equation}
see \cite{He16}.
Then, we denote by
\begin{equation*}
    \mathbb{NN}^{\sigma,l,m}_{d_1, d_2} :=
    \bigl\{
        L_l \circ \dots \circ L_0 : \real^{d_1} \to \real^{d_2}
    \ : \
        L_0 \in \mathbb{L}^{\sigma}_{d_1, m},
        L_l \in \mathbb{L}^{\sigma_{\rm Id}}_{m, d_2},
        L_i \in \mathbb{L}^{\sigma, \rm res}_m,
        1 \leq i < l
    \bigr\}
\end{equation*}
the set of feed-forward neural networks
with one output layer,
$l \geq 1$ hidden residual layers
each containing $m \geq 1$ neurons,
and the activation functions of
the output layer and
the hidden layers
being respectively
the identity function $\sigma_{\rm Id}$
and $\sigma$.
Any $v(\cdot; \theta) \in \mathbb{NN}^{\sigma,l,m}_{d_1, d_2}$
is fully determined by the sequence
\begin{equation*}
    \theta := \bigl( W_0, b_0, W_1, b_1, \dots, W_{l-1}, b_{l-1}, W_l, b_l \bigr)
\end{equation*}
of $\left( (d_1+1) m + (l - 1) (m+1) m + (m+1) d_2 \right)$ parameters.

\medskip

Since by the universal approximation theorem,
see e.g.\ Theorem~1 of \cite{hornik1991approximation},
$\bigcup\limits_{m = 1}^{\infty} \mathbb{NN}^{\sigma,l ,m}_{d_1, d_2}$
is dense in $L^2$ functional space,
the optimization problem \eqref{eq:L2 minimality} can be approximated by
\begin{equation}
\label{eq:NN approximation}
    v^* \approx \argmin\limits_{v \in \mathbb{NN}^{\sigma,l ,m}_{d + 1, d}}
        \sum\limits_{i = 1}^d
        \E \left[
            \left(\mathcal{H}\left({\zeta , X, {\rm Id}_i }\right)
                    - v_i(\zeta , X)\right)^2 \right].
\end{equation}
By the law of large numbers,
\eqref{eq:NN approximation} can be further approximated by
\begin{equation}
\label{eq:monte carlo}
    v^* \approx \argmin\limits_{v \in \mathbb{NN}^{\sigma,l ,m}_{d + 1, d}}
        \sum\limits_{i = 1}^d
        N^{-1} \sum\limits_{j = 1}^N
        \left(\mathcal{H}_{i,j} - v_i(\zeta_j, X_j)\right)^2,
\end{equation}
where for all $j = 1, \dots, N$,
$(\zeta_j, X_j)$ is drawn independently from the distribution of $(\zeta , X)$
and $\mathcal{H}_{i,j}$ is drawn from
$\mathcal{H}_{\zeta_j, X_j, {\rm Id}_i }$
using Algorithm~\ref{alg:coding tree}.
However,
the approximation \eqref{eq:monte carlo} may perform poorly
when the variance of $\mathcal{H}_{i,j}$ is too high.
To address this issue,
we perform
\begin{equation}
\label{eq:monte carlo more samples}
    v^* \approx \argmin\limits_{v \in \mathbb{NN}^{\sigma,l ,m}_{d + 1, d}}
        \sum\limits_{i = 1}^d
        \frac{1}{N} \sum\limits_{j = 1}^N
            \left(\frac{1}{M} \sum\limits_{k = 1}^M \mathcal{H}_{i,j,k}
                    - v_i(\zeta_j, X_j)\right)^2,
\end{equation}
where for all $k = 1, \dots, M$,
$\mathcal{H}_{i,j,k}$ is drawn independently
from $\mathcal{H}_{\zeta_j, X_j, {\rm Id}_i }$
using Algorithm~\ref{alg:coding tree}.

\medskip

Finally, the deep branching algorithm
using the gradient descent method
to solve the optimization in \eqref{eq:monte carlo more samples}
is summarized in Algorithm~\ref{alg:deep}.
      
\begin{algorithm}[H]
  \caption{Deep branching algorithm.}
  \label{alg:deep}
    \begin{algorithmic}
        \Require The learning rate $\eta$ and the number of epochs $P$
        \Ensure $v(\cdot, \cdot; \theta) \in \mathbb{NN}^{\sigma,l,m}_{d+1, d}$
        \State $(\zeta_j, X_j)_{1 \leq j \leq N} \gets$ random vectors
                    drawn from the distribution of $(\zeta , X)$
        \State $(\mathcal{H}_{i,j,k})_{\footnotesize
                    \substack{1 \leq i \leq d \\
                              1 \leq j \leq N \\
                              1 \leq k \leq M}} \gets$
                    random variables
                    generated by TREE$(\zeta_j, X_j, {\rm Id}_i )$ in
                    Algorithm~\ref{alg:coding tree}
        \smallskip
        \State Initialize $\theta$
        \For{$i \gets 1, \dots, P$}
            \State $L \gets
                \sum\limits_{i = 1}^d
                N^{-1} \sum\limits_{j = 1}^N
                \bigg( M^{-1} \sum\limits_{k = 1}^M \mathcal{H}_{i,j,k}
                    - v_i(\zeta_j, X_j; \theta)\bigg)^2$
            \State $\theta \gets \theta - \eta \nabla_\theta L$
        \EndFor
    \end{algorithmic}
\end{algorithm}
\noindent
     Since no closed form expression may be available for the function
          $$
          \phi_0(x)
          =
          \frac{\Gamma (d/2)}{2 \pi^{d/2}}
          \bigintsss_{\real^d}
          \frac{N(y)}{\abs{y}^d}
          f_0\big(
                   \partial_{\bar{\alpha}^{q+1}}\phi_{\beta_{q+1}}(x+y)
                 ,
                 \ldots ,
                            \partial_{\bar{\alpha}^n}\phi_{\beta_n}(x+y)
                 \big)
             dy,
          $$
          we approximate it using the neural network function
          and Monte Carlo method for the numerical integration.
          More precisely, we approximate
          \begin{equation}
          \label{eq:phi0 mc}
          \phi_0(x) =
          \int_0^\infty
          \int_{\real^d}
          \varphi_1 (s, y)
  \frac{N(y)}{2s}
          f_0\big(
          \partial_{\bar{\alpha}^{q+1}}\phi_{\beta_{q+1}}(x+y), \ldots ,
          \partial_{\bar{\alpha}^n}\phi_{\beta_n}(x+y)
          \big)
          dy
          ds
          \end{equation}
          using
          \begin{align*}
          \phi_0 \approx
          \argmin\limits_{v \in \mathbb{NN}^{\sigma,l ,m}_{d, 1}}
          N^{-1} \sum\limits_{i = 1}^N \biggl(
          &  M^{-1} \sum\limits_{j = 1}^M
                  N(Y_{i,j}) \times
                  (2 \widetilde{\tau}_{i,j} \ \!
                  \widetilde{\rho}(\widetilde{\tau_{i,j}}))^{-1}
                  \times
          \\
          & f_0\big(
              \partial_{\bar{\alpha}^{q+1}}\phi_{\beta_{q+1}}(X_i+Y_{i,j})
                 ,
                 \ldots ,
                 \partial_{\bar{\alpha}^n}\phi_{\beta_n}(X_i+Y_{i,j})
                 \big)
              - v_i(\zeta_i, X_i)\biggr)^2,
          \end{align*}
          where
          $X_i$ is the uniform vector on
          $[x_{\rm min} - ( x_{\rm max} - x_{\rm min} ) / 2,
          x_{\rm max} + ( x_{\rm max} - x_{\rm min} ) / 2 ]^d$,
          $\widetilde{\tau}_{i, j}$ is the random variable
          drawn independently from the distribution of $\widetilde{\rho}$,
          and $Y_{i, j}$ is the random vector
          drawn independently from $\mathcal{N}(0, \widetilde{\tau}_{i, j})$,
          see \eqref{eq:poisson integration calculation}
          for the derivation of \eqref{eq:phi0 mc}.

          \medskip

          \noindent
           Algorithm~\ref{alg:deep} is implemented with the following parameters:
   \begin{enumerate}[label=\emph{\alph*})]
         \item $\rho$ is chosen to be the PDF of exponential distribution
                with rate $- ( \log 0.95) / T$,
        \item $\widetilde{\rho}$ is chosen to be the PDF of uniform distribution
          $\widetilde{\rho}(x) =
          (6 - 10^{-5})^{-1}
          \bm{1}_{[10^{-5}, 6]}(x)
          $,
        \item given $x_{\rm min} < x_{\rm max}$,
 we let $(\zeta , X)$ be a uniformly distributed random vector on
 $[0, T] \times \Omega$,
 where $    \Omega := [x_{\rm min}, x_{\rm max}]^d$,
\item
 the activation function 
 $\sigma_{\tanh}(x) := \tanh(x)$ is used
 instead of ReLu 
 because the target PDE solution \eqref{eq:main pde} is smooth,
\item the optimal learning rate $\eta$ for gradient update 
is obtained by trial and error,
given that that a lower $\eta$
      means slow convergence to a possibly local suboptimum,
      while a higher $\eta$ can lead to instability,
    \item
      standard parameters without tuning were used for Adam optimization
      and batch normalization,
   \end{enumerate}
 and we perform the following additional steps:
    \begin{enumerate}[label=\emph{\alph*}),resume]
        \item $\eta \gets \eta / 10$ at epoch $1,000$ and $2,000$.
        \item Instead of using $\eta$ to update $\theta$ directly,
                the Adam algorithm is used to update $\theta$,
                see \cite{kingma2014adam}.
        \item A batch normalization layer
                is added before the every layer of
                                \eqref{eq:one layer}-\eqref{eq:one layer resnet},
                                see \cite{ioffe2015batch}.
\end{enumerate}

\section{Application to the Navier-Stokes equation}
            \label{sec:numerical examples}
 The incompressible Navier-Stokes equation 
\begin{equation}
\nonumber 
\begin{cases}
  \displaystyle
  \partial_t u_i(t,x) + \nu \Delta u_i(t,x)
  =
  \partial_{\bm{1}_i} p(t,x)
 + \sum\limits_{j=1}^d u_j(t,x) \partial_{\bm{1}_j} u_i(t,x),
  \medskip
  \\
  \displaystyle
  \Delta p(t, x)
  = -\sum\limits_{i,j=1}^d
\partial_{\bm{1}_j} u_i(t,x)
\partial_{\bm{1}_i} u_j(t,x),
  \medskip
  \\
u_i(T,x) = \phi_i (x),
\quad (t,x) = (t,x_1, \ldots, x_d) \in [0,T] \times \real^d,
\quad
 i = 1,\ldots , d,
\end{cases}
\end{equation}
with pressure term $p(t,x)=u_0(t,x)$
can be obtained as a particular case of the system \eqref{eq:main pde}.
For this, we take $n=d(d+2)$, $q=d$, and let
$$
f_0 \big(
y_1,\ldots , y_d,
z^{(1)}_1,\ldots , z^{(1)}_d,
\ldots ,
z^{(d)}_1,\ldots , z^{(d)}_d
\big)
= -\sum\limits_{i,j=1}^d
z_i^{(j)}
z_j^{(i)}
$$
and
$$
 f_i\big(
 x_1,\ldots , x_d,
 y_1,\ldots , y_d,
 z^{(1)}_1,\ldots , z^{(1)}_d,
 \ldots ,
 z^{(d)}_1,\ldots , z^{(d)}_d
 \big)
 = - x_i - \sum\limits_{j=1}^d y_j z_i^{(j)},
$$
 $i = 1 , \ldots , d$, with
 $\bar{\alpha}^i = \bm{1}_i$, $i=1,\ldots , d$,
 $\bar{\alpha}^{d+1} = \cdots = \bar{\alpha}^{2d} = 0$,
 $\beta_{d+i} = i$, $i=1,\ldots ,d$,
 $\bar{\alpha}^{i+(j+1)d} = \bm{1}_j$,
 $\beta_{i+(j+1)d} = i$, $i,j=1,\ldots , d$, so that
\begin{eqnarray*}
 f_0\big(
          \partial_{\bar{\alpha}^{q+1}}u_{\beta_{q+1}}(t,x)
        ,
        \ldots ,
                   \partial_{\bar{\alpha}^n}u_{\beta_n}(t,x)
        \big)
 & = &
        f_0\big(
                u_1(t,x),\ldots , u_d(t,x) ,
(\partial_{\bm{1}_j}u_k(t,x))_{
  1 \leq j,k \leq d}
\big)
\\
 & = & -\sum\limits_{i,j=1}^d
\partial_{\bm{1}_j} u_i(t,x)
\partial_{\bm{1}_i} u_j(t,x),
\end{eqnarray*}
and
\begin{eqnarray*}
\lefteqn{ f_i\big(
    \partial_{\bar{\alpha}^1}u_0 (t,x)
    ,
    \ldots ,
    \partial_{\bar{\alpha}^q }u_0 (t,x)
    ,
        \partial_{\bar{\alpha}^{q+1}}u_{\beta_{q+1}}(t,x)
        ,
        \ldots ,
                   \partial_{\bar{\alpha}^n}u_{\beta_n}(t,x)
                   \big)
}
\\
 & = &
f_i\big(
\partial_{\bm{1}_1} u_0(t,x)
,
\ldots
,
\partial_{\bm{1}_d} u_0(t,x)
,
u_1(t,x),\ldots , u_d(t,x)
,
(\partial_{\bm{1}_j}u_k(t,x))_{
                1 \leq j,k \leq d}
\big)
\\
 & = &
-\partial_{\bm{1}_i} u_0(t,x)
-\sum\limits_{j=1}^d u_j(t,x) \partial_{\bm{1}_j} u_i(t,x).
\end{eqnarray*}
 The following numerical examples in Sections~\ref{sub1}-\ref{sub2}
 are implemented in {\sc Python} using {\sc PyTorch}
 on a computer with a 3.60 GHz AMD Ryzen 5 3500 
  processor, a 16 GB at 3200 MHz
 DDR4-SDRAM, and a GeForce RTX 3080 Ti graphics card with 12 GB memory.
The default {\sc PyTorch} initialization scheme for $\theta$ is used,
together with the default values
$N=100,000$, $M=1,000$, $P=10,000$, $\eta=0.01$, $l=3$,
$m=100$, $x_{\rm min}=0$, $x_{\rm max}=2\pi$.
For any $\delta > 0$,
we let $C_\delta := \delta \zz^d$
and perform the analysis of error on the grid of
$\Omega \cap C_\delta$ at time $t_k = {kT} / {10}$
for $k = 0, 1, \ldots, 9$.
Our benchmarking to \cite{angeli} and \cite{lejay2020forward}
 uses the following errors:
\begin{align*}
& e_i(t_k)
= \sup\limits_{x \ \! \in \ \! \Omega \ \! \cap \ \! C_\delta}
\abs{u_i(t_k, x) - v_i(t_k, x; \theta)}^2,
\\
& e(t_k)
= \sup\limits_{x \ \! \in \ \! \Omega \ \! \cap \ \! C_\delta}
\sum\limits_{i = 1}^d \abs{u_i(t_k, x) - v_i(t_k, x; \theta)}^2,
\\
& {\rm erru}(t_k)
=
\left(
\frac
{\sum\limits_{i = 1}^d
\sum\limits_{x \ \! \in \ \! \Omega \ \! \cap \ \! C_\delta}
\abs{u_i(t_k, x) - v_i(t_k, x; \theta)}^2}
{\sum\limits_{i = 1}^d
\sum\limits_{x \ \! \in \ \! \Omega \ \! \cap \ \! C_\delta}
\abs{u_i(t_k, x)}^2}
\right)^{1/2},
\\
& {\rm errgu}(t_k)
= \left(
\frac
{\sum\limits_{i, j = 1}^d
\sum\limits_{x \ \! \in \ \! \Omega \ \! \cap \ \! C_\delta}
\abs{\partial_{\bm{1}_j} u_i(t_k, x) - \partial_{\bm{1}_j}v_i(t_k, x; \theta)}^2}
{\sum\limits_{i, j = 1}^d
\sum\limits_{x \ \! \in \ \! \Omega \ \! \cap \ \! C_\delta}
\abs{\partial_{\bm{1}_j} u_i(t_k, x)}^2}
\right)^{1/2},
\\
& {\rm errdivu}(t_k)
=
\left(
(x_{\rm max} - x_{\rm min})^d
\abs{\Omega}^{-1}
\sum\limits_{i = 1}^d
\sum\limits_{x \ \! \in \ \! \Omega \ \! \cap \ \! C_\delta}
\abs{\partial_{\bm{1}_i} u_i(t_k, x)}^2
\right)^{1/2},
\\
& {\rm errp}(T)
=
\left(
\frac
{\sum\limits_{x \ \! \in \ \! \Omega \ \! \cap \ \! C_\delta}
\abs{p(T, x)
- v_0(T, x; \theta)
+ \abs{\Omega}^{-1}
\sum\limits_{x \ \! \in \ \! \Omega \ \! \cap \ \! C_\delta} v_0(T, x; \theta)}^2}
{\sum\limits_{x \ \! \in \ \! \Omega \ \! \cap \ \! C_\delta}
\abs{p(T, x)}^2}
\right)^{1/2}
.
\end{align*}

\subsection{Taylor-Green vortex}
\label{sub1}
\begin{figure}[H]
  \centering
 \begin{subfigure}[b]{0.45\textwidth}
    \includegraphics[width=\linewidth]{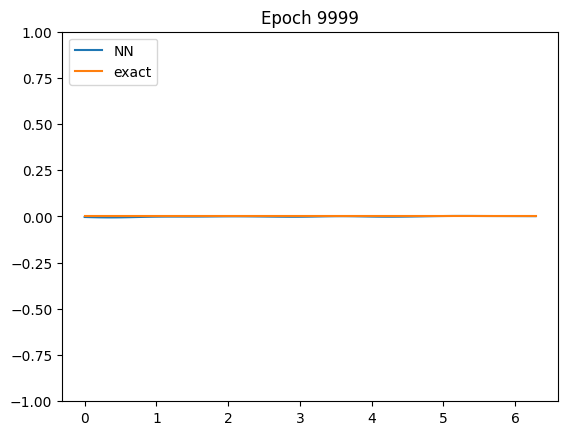}
    \caption{Comparison for $u_1 (x_1,\pi)$.}
 \end{subfigure}
  \begin{subfigure}[b]{0.45\textwidth}
    \includegraphics[width=\linewidth]{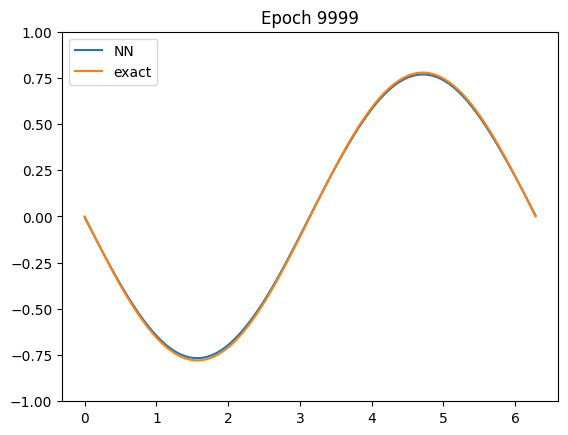}
    \caption{Comparison for $u_2 (x_1,\pi)$.}
 \end{subfigure}
   \vskip-0.1cm
\caption{Comparison with the exact solution \eqref{eq:taylor green} by taking $t = {T}/{2}$ and $x_2 = \pi$.}
\label{fig:taylor green 1}
\end{figure}

In this section we consider the $2$-dimensional Taylor-Green \cite{taylor1937mechanism} vortex
\begin{equation}
\label{eq:taylor green}
\begin{cases}
  \displaystyle
  u_1(t, x) = - \cos(x_1) \sin(x_2) \re^{-2\nu (T - t)},
  \medskip
  \\
  u_2(t, x) = \sin(x_1) \cos(x_2) \re^{-2\nu (T - t)},
  \medskip
  \\
  \displaystyle
  u_0(t, x) = -\frac{1}{4} \left(\cos(2x_1) + \cos(2x_2)\right) \re^{-4\nu (T - t)} + c
\end{cases}
\end{equation}
$x=(x_1,x_2) \in [0,2\pi]^2$,
with Reynolds numbers in the range $[1,100]$.
 We first let $\nu = 1$, $\delta= \pi / 126$ , $T=1/4$, and present the results
in Figure~\ref{fig:taylor green 1} and Table~\ref{table:taylor green 1}.
 In this example and the next one, our method provides a
       solution on $[0,T]\times \real^d$
       by only imposing a terminal condition at terminal time
       $t=T$. 
       As those examples are periodic we only
       provide solution values on a given interval of periodicity
       as in \cite{lejay2020forward},
       however our solver can be used to yield estimates on larger intervals as
       well. 

\vspace{-0.4cm}

\begin{table}[H]
    \centering
    \resizebox{\textwidth}{!}{\begin{tabular}{|r|c|c|c|c|c|c|c|c|c|c|c|}
        \hline
\multicolumn{1}{|c|}{k}
        & 0 & 1 & 2 & 3 & 4 & 5 & 6 & 7 & 8 & 9 & 10\\
        \hline
${e}_0(t_k)$
& 1.90E-04 & 1.94E-04 & 1.91E-04 & 2.03E-04 & 2.51E-04 & 2.96E-04 & 3.28E-04 & 3.53E-04 & 3.79E-04 & 4.64E-04 & --- \\
${e}_1(t_k)$
& 1.99E-04 & 1.95E-04 & 2.24E-04 & 2.34E-04 & 2.44E-04 & 2.56E-04 & 2.63E-04 & 2.44E-04 & 2.40E-04 & 3.63E-04 & --- \\
${e}(t_k)$
& 2.32E-04 & 1.98E-04 & 2.27E-04 & 2.40E-04 & 2.58E-04 & 3.06E-04 & 3.39E-04 & 3.64E-04 & 3.86E-04 & 4.65E-04 & --- \\
\hline
erru($t_k$)
& 1.57E-02 & 1.43E-02 & 1.51E-02 & 1.64E-02 & 1.73E-02 & 1.76E-02 & 1.72E-02 & 1.64E-02 & 1.56E-02 & 1.56E-02 & --- \\
errgu($t_k$)
& 3.24E-02 & 2.75E-02 & 2.50E-02 & 2.40E-02 & 2.37E-02 & 2.34E-02 & 2.28E-02 & 2.20E-02 & 2.14E-02 & 2.20E-02 & --- \\
errdivu($t_k$)
& 2.03E-02 & 1.38E-02 & 9.92E-03 & 9.13E-03 & 1.02E-02 & 1.12E-02 & 1.14E-02 & 1.06E-02 & 9.81E-03 & 1.33E-02 & --- \\
errp($t_k$)
& --- & --- & --- & --- & --- & --- & --- & --- & --- & --- & 1.75E-02 \\
        \hline
\end{tabular}}
		\caption{Error comparison.}
        \label{table:taylor green 1}
\end{table}

\vspace{-0.4cm}

\noindent
 Our simulation runtime on the full grid $[0, T] \times \Omega$
 is approximately 22 minutes for the Taylor-Green vortex, 
 after 20 minutes of pre-computation for the training of
 the terminal condition $p(T,x)$. 
 Table~\ref{table:taylor green 1} above can be compared\footnote{The numbers in Table~\ref{table:taylor green 1} above should be multiplied by $10^3$ for comparison with Table~1 in \cite{lejay2020forward}.} to Table~1 in \cite{lejay2020forward} where computing a single time step by BSDEs and Monte Carlo 
 on a computer cluster with a few tens of cores
 took approximately 2 hours, whereas our neural network approach yields a functional estimate on $[0,2\pi]^2 \times [0,T]$.

\medskip

\begin{figure}[H]
  \centering
 \begin{subfigure}[b]{0.45\textwidth}
    \includegraphics[width=\linewidth]{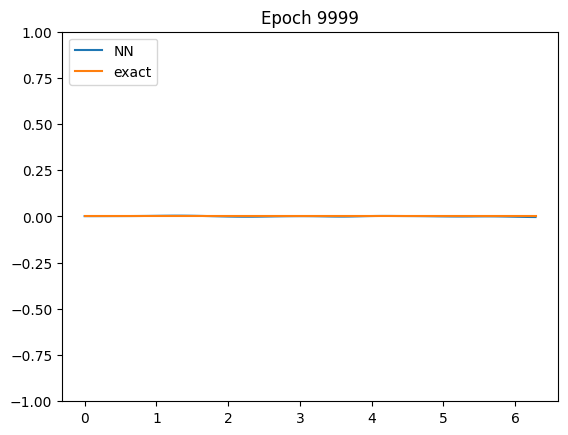}
    \caption{Comparison for $u_1 (x_1,\pi)$.}
 \end{subfigure}
  \begin{subfigure}[b]{0.45\textwidth}
    \includegraphics[width=\linewidth]{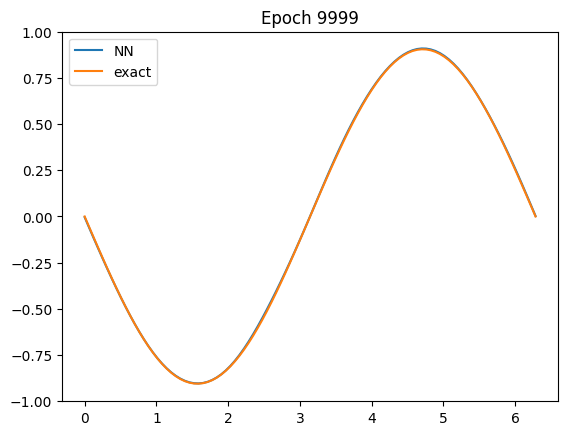}
    \caption{Comparison for $u_2 (x_1,\pi)$.}
 \end{subfigure}
   \vskip-0.1cm
   \caption{Comparison with the exact solution \eqref{eq:taylor green} by taking $t = {T}/{2}$ and $x_2 = \pi$.}
\label{fig:taylor green 2}
\end{figure}

Next, we let $\nu = 0.1$, $T=1$, and present the results
in Figure~\ref{fig:taylor green 2} and Table~\ref{table:taylor green 2}.

\vskip-0.4cm

\begin{table}[H]
    \centering
    \resizebox{\textwidth}{!}{\begin{tabular}{|r|c|c|c|c|c|c|c|c|c|c|c|}
        \hline
\multicolumn{1}{|c|}{k}
        & 0 & 1 & 2 & 3 & 4 & 5 & 6 & 7 & 8 & 9 & 10\\
        \hline
${e}_0(t_k)$
& 3.27E-04 & 2.48E-04 & 1.73E-04 & 1.41E-04 & 1.53E-04 & 1.62E-04 & 1.67E-04 & 1.73E-04 & 1.89E-04 & 2.07E-04 & --- \\
${e}_1(t_k)$
& 3.27E-04 & 2.04E-04 & 1.34E-04 & 1.28E-04 & 1.17E-04 & 1.12E-04 & 1.26E-04 & 1.54E-04 & 1.93E-04 & 2.54E-04 & --- \\
${e}(t_k)$
& 3.84E-04 & 2.72E-04 & 1.79E-04 & 1.67E-04 & 1.75E-04 & 1.78E-04 & 1.75E-04 & 1.77E-04 & 2.14E-04 & 2.70E-04 & --- \\
\hline
erru($t_k$)
& 1.29E-02 & 1.04E-02 & 8.90E-03 & 8.34E-03 & 8.27E-03 & 8.28E-03 & 8.18E-03 & 8.03E-03 & 8.17E-03 & 9.18E-03 & --- \\
errgu($t_k$)
& 3.03E-02 & 2.56E-02 & 2.17E-02 & 1.88E-02 & 1.70E-02 & 1.61E-02 & 1.62E-02 & 1.73E-02 & 1.94E-02 & 2.25E-02 & --- \\
errdivu($t_k$)
& 2.52E-02 & 1.76E-02 & 1.32E-02 & 1.21E-02 & 1.29E-02 & 1.35E-02 & 1.34E-02 & 1.29E-02 & 1.36E-02 & 1.78E-02 & --- \\
errp($t_k$)
& --- & --- & --- & --- & --- & --- & --- & --- & --- & --- & 1.75E-02 \\
        \hline
\end{tabular}}
		\caption{Error comparison.}
        \label{table:taylor green 2}
\end{table}

\vspace{-0.4cm}

\noindent
Table~\ref{table:taylor green 2} above can be compared to Tables~17 and 19 in Section~5 of \cite{angeli}, which use mesh-based methods running a 20 core CPU under Ubuntu 16.04 with 32 Go RAM.
Our results are comparable in terms of errgu$(t_k)$ to the rectangular meshes 1 to 4 in Table~19 therein, which require up to 5 seconds. 
Those results are also comparable in terms of erru$(t_k)$ to the triangular meshes 1 to 3 in Table~17 therein, which require up to 44 seconds. 

\medskip

\begin{figure}[H]
  \centering
 \begin{subfigure}[b]{0.45\textwidth}
    \includegraphics[width=\linewidth]{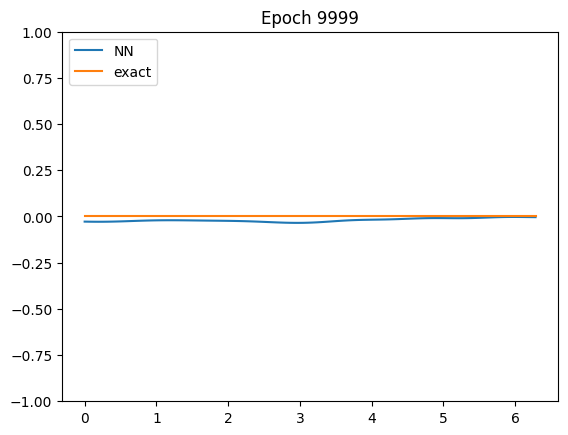}
    \caption{Comparison for $u_1 (x_1,\pi)$.}
 \end{subfigure}
  \begin{subfigure}[b]{0.45\textwidth}
    \includegraphics[width=\linewidth]{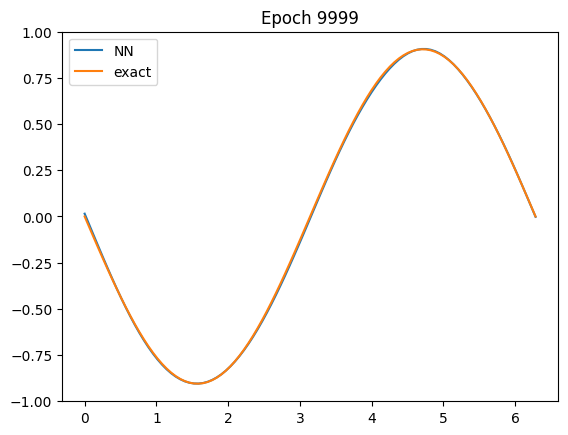}
    \caption{Comparison for $u_2 (x_1,\pi)$.}
 \end{subfigure}
   \vskip-0.1cm
   \caption{Comparison with the exact solution \eqref{eq:taylor green} by taking $t = {T}/{2}$ and $x_2 = \pi$.}
   \label{fig:taylor green 3}
\end{figure}

\noindent
We now let $\nu = 0.01$, $T=10$, and present the results
in Figure~\ref{fig:taylor green 3} and Table~\ref{table:taylor green 3}.

\vskip-0.4cm

\begin{table}[H]
    \centering
    \resizebox{\textwidth}{!}{\begin{tabular}{|r|c|c|c|c|c|c|c|c|c|c|c|}
        \hline
\multicolumn{1}{|c|}{k}
        & 0 & 1 & 2 & 3 & 4 & 5 & 6 & 7 & 8 & 9 & 10\\
        \hline
${e}_0(t_k)$
& 4.99E-03 & 4.51E-03 & 4.01E-03 & 3.20E-03 & 3.16E-03 & 2.28E-03 & 1.74E-03 & 1.33E-03 & 1.56E-03 & 1.69E-03 & --- \\
${e}_1(t_k)$
& 2.54E-03 & 2.14E-03 & 1.73E-03 & 1.27E-03 & 8.45E-04 & 5.77E-04 & 4.96E-04 & 5.29E-04 & 8.22E-04 & 1.63E-03 & --- \\
${e}(t_k)$
& 5.55E-03 & 4.94E-03 & 4.29E-03 & 3.39E-03 & 3.39E-03 & 2.42E-03 & 1.84E-03 & 1.42E-03 & 1.56E-03 & 1.88E-03 & --- \\
\hline
erru($t_k$)
& 6.22E-02 & 5.55E-02 & 4.94E-02 & 4.36E-02 & 3.83E-02 & 3.36E-02 & 2.97E-02 & 2.72E-02 & 2.65E-02 & 2.81E-02 & --- \\
errgu($t_k$)
& 4.39E-02 & 3.86E-02 & 3.43E-02 & 3.09E-02 & 2.85E-02 & 2.72E-02 & 2.72E-02 & 2.87E-02 & 3.16E-02 & 3.59E-02 & --- \\
errdivu($t_k$)
& 6.28E-02 & 5.11E-02 & 4.20E-02 & 3.48E-02 & 2.90E-02 & 2.48E-02 & 2.33E-02 & 2.59E-02 & 3.23E-02 & 4.16E-02 & --- \\
errp($t_k$)
& --- & --- & --- & --- & --- & --- & --- & --- & --- & --- & 1.75E-02 \\
        \hline
\end{tabular}}
		\caption{Error comparison.}
        \label{table:taylor green 3}
\end{table}

\vspace{-0.4cm}

\noindent
Table~\ref{table:taylor green 3} above is comparable in terms of errgu$(t_k)$ to the rectangular meshes 1 to 4 in Table~20 in Section~5 of \cite{angeli}, which require up to 10 seconds. 
Those results are also comparable in terms of erru$(t_k)$ to the triangular meshes 1 to 3 in Table~18 therein, which require up to one minute. 

\subsection{Arnold-Beltrami-Childress flow}
\label{sub2}
\begin{figure}[H]
  \centering
 \begin{subfigure}[b]{0.3\textwidth}
    \includegraphics[width=\linewidth]{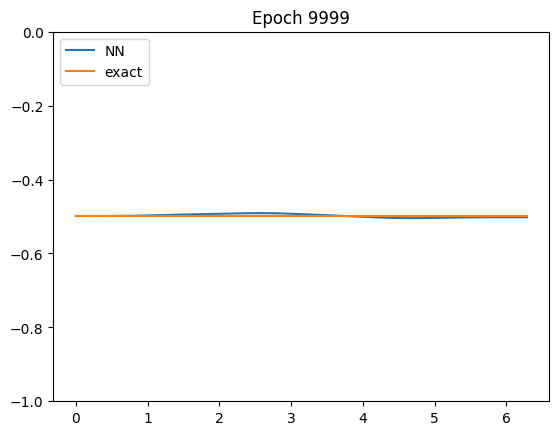}
    \caption{Comp. for $u_1 (x_1,\pi,\pi)$.}
 \end{subfigure}
  \begin{subfigure}[b]{0.3\textwidth}
    \includegraphics[width=\linewidth]{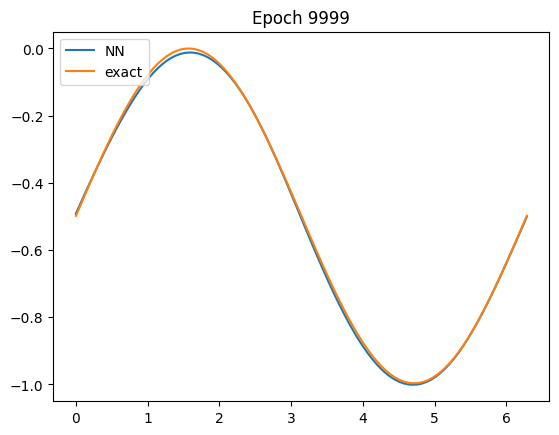}
    \caption{Comp. for $u_2 (x_1,\pi,\pi)$.}
 \end{subfigure}
 \begin{subfigure}[b]{0.3\textwidth}
    \includegraphics[width=\linewidth]{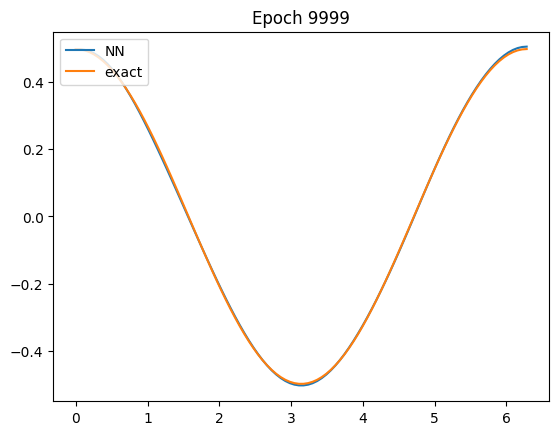}
    \caption{Comp. for $u_3 (x_1,\pi,\pi)$.}
 \end{subfigure}
   \vskip-0.1cm
   \caption{Comparison with the exact solution \eqref{eq:abc flow} by taking $t = {T}/{2}$ and $x_2 = x_3 = \pi$.}
\label{fig:abc}
\end{figure}

Here, we consider the following
$3$-dimensional Arnold-Beltrami-Childress \cite{arnold1965topologie},
 \cite{childress1970new} flow
\begin{equation}
\label{eq:abc flow}
\begin{cases}
  \displaystyle
  u_1(t, x) = \left(A \sin(x_3) + C \cos(x_2)\right) \re^{-\nu (T - t)},
  \medskip
  \\
  u_2(t, x) = \left(B \sin(x_1) + A \cos(x_3)\right) \re^{-\nu (T - t)},
  \medskip
  \\
  u_3(t, x) = \left(C \sin(x_2) + B \cos(x_1)\right) \re^{-\nu (T - t)},
  \medskip
  \\
  u_0(t, x) = -\left(AC \sin(x_3) \cos(x_2) + BA \sin(x_1) \cos(x_3) + CB \sin(x_2) \cos(x_1)\right) \re^{-2\nu (T - t)} + c.
\end{cases}
\end{equation}
$x=(x_1,x_2,x_3) \in [0,2\pi]^3$.
We first let $\nu = 0.01$, $A=B=C=0.5$, $T=0.7$, $\delta= \pi / 45$,
which corresponds to a Reynolds numbers in the range $[1,100]$,
and we present the results
in Figure~\ref{fig:abc} and Table~\ref{table:abc}.

\vskip-0.4cm

\begin{table}[H]
    \centering
    \resizebox{\textwidth}{!}{\begin{tabular}{|r|c|c|c|c|c|c|c|c|c|c|c|}
        \hline
\multicolumn{1}{|c|}{k}
        & 0 & 1 & 2 & 3 & 4 & 5 & 6 & 7 & 8 & 9 & 10 \\
        \hline
${e}_0(t_k)$
& 9.68E-04 & 4.27E-04 & 5.92E-04 & 5.29E-04 & 6.02E-04 & 5.42E-04 & 5.61E-04 & 4.11E-04 & 6.28E-04 & 5.96E-04 & --- \\
${e}_1(t_k)$
& 2.20E-03 & 1.01E-03 & 8.25E-04 & 9.41E-04 & 7.28E-04 & 7.18E-04 & 7.90E-04 & 8.11E-04 & 6.64E-04 & 7.05E-04 & --- \\
${e}_2(t_k)$
& 1.09E-03 & 5.06E-04 & 4.61E-04 & 6.83E-04 & 9.28E-04 & 4.70E-04 & 6.65E-04 & 8.29E-04 & 4.90E-04 & 5.68E-04 & --- \\
${e}(t_k)$
& 2.97E-03 & 1.34E-03 & 1.12E-03 & 1.50E-03 & 1.20E-03 & 1.05E-03 & 1.04E-03 & 1.10E-03 & 6.99E-04 & 9.08E-04 & --- \\
\hline
erru($t_k$)
& 1.64E-02 & 1.14E-02 & 1.18E-02 & 1.22E-02 & 1.32E-02 & 1.28E-02 & 1.24E-02 & 1.14E-02 & 1.11E-02 & 1.09E-02 & --- \\
errgu($t_k$)
& 3.72E-02 & 3.16E-02 & 3.19E-02 & 3.34E-02 & 3.59E-02 & 3.52E-02 & 3.41E-02 & 3.15E-02 & 3.01E-02 & 3.02E-02 & --- \\
errdivu($t_k$)
& 8.74E-02 & 5.42E-02 & 5.30E-02 & 5.84E-02 & 6.91E-02 & 7.12E-02 & 6.02E-02 & 5.45E-02 & 5.15E-02 & 5.05E-02 & --- \\
errp($t_k$)
& --- & --- & --- & --- & --- & --- & --- & --- & --- & --- & 1.93E-02 \\
        \hline
\end{tabular}}
		\caption{Error comparison.}
        \label{table:abc}
\end{table}

\vspace{-0.4cm}

\noindent
Our simulation runtime on the full grid $[0, T] \times \Omega$
is approximately 60 minutes for the Arnold-Beltrami-Childress flow
after 30 minutes of pre-computation for the training of
the terminal condition $p(T,x)$. 
Table~\ref{table:abc} can be compared\footnote{The numbers in Table~\ref{table:abc} above should be multiplied by $10^2$ for comparison with Table~5 in \cite{lejay2020forward}.} to Table~5 in \cite{lejay2020forward} where a single time step
by BSDEs and Monte Carlo took approximately 20 hours.
Our results have a significantly lower runtime,
and are at least one order of magnitude more accurate
than \cite{lejay2020forward}.
In addition, the neural network approach yields a functional estimate on
$[0,2\pi]^3 \times [0,T]$ instead of estimating the solution at discrete time
instants.

\medskip

\begin{figure}[H]
  \centering
 \begin{subfigure}[b]{0.3\textwidth}
    \includegraphics[width=\linewidth]{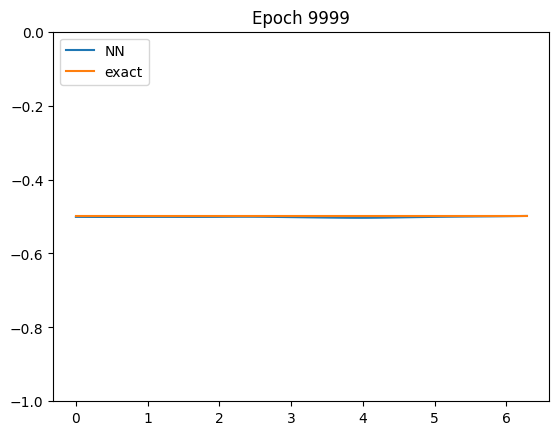}
    \caption{Comp.for $u_1(x_1,\pi,\pi)$.}
 \end{subfigure}
  \begin{subfigure}[b]{0.3\textwidth}
    \includegraphics[width=\linewidth]{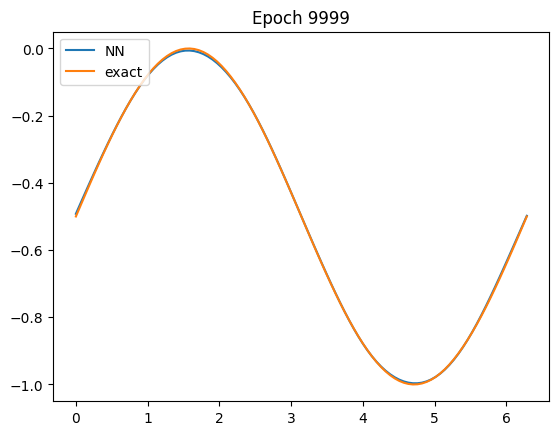}
    \caption{Comp. for $u_2(x_1,\pi,\pi)$.}
 \end{subfigure}
 \begin{subfigure}[b]{0.3\textwidth}
    \includegraphics[width=\linewidth]{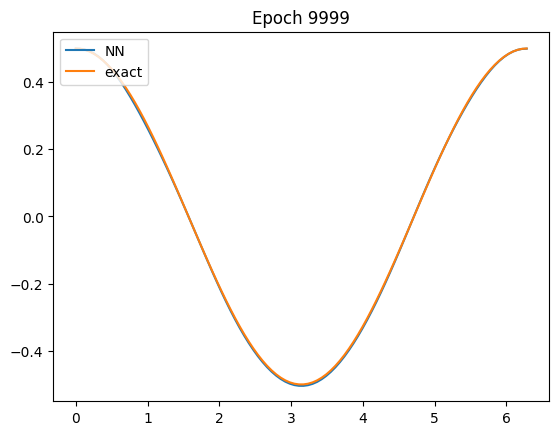}
    \caption{Comp. for $u_3(x_1,\pi,\pi)$.}
 \end{subfigure}
   \vskip-0.1cm
   \caption{Comparison with the exact solution \eqref{eq:abc flow} by taking $t = {T}/{2}$ and $x_2 = x_3 = \pi$.}
\label{fig:abc low beta}
\end{figure}

Finally, we let $\nu=10^{-4}$, which corresponds to a Reynolds number of order 10,000,
and we present the results in Figure~\ref{fig:abc low beta} and Table~\ref{table:abc low beta}.

\vskip-0.4cm

\begin{table}[H]
    \centering
    \resizebox{\textwidth}{!}{\begin{tabular}{|r|c|c|c|c|c|c|c|c|c|c|c|}
        \hline
\multicolumn{1}{|c|}{k}
        & 0 & 1 & 2 & 3 & 4 & 5 & 6 & 7 & 8 & 9 & 10 \\
        \hline
${e}_0(t_k)$
& 2.86E-04 & 1.46E-04 & 1.29E-04 & 1.15E-04 & 1.75E-04 & 1.47E-04 & 1.73E-04 & 9.17E-05 & 1.31E-04 & 1.15E-04 & --- \\
${e}_1(t_k)$
& 3.22E-04 & 1.58E-04 & 1.56E-04 & 1.67E-04 & 1.27E-04 & 1.37E-04 & 1.14E-04 & 1.11E-04 & 1.55E-04 & 2.12E-04 & --- \\
${e}_2(t_k)$
& 2.68E-04 & 1.72E-04 & 2.48E-04 & 2.09E-04 & 1.72E-04 & 1.73E-04 & 1.71E-04 & 1.42E-04 & 1.15E-04 & 1.52E-04 & --- \\
${e}(t_k)$
& 4.05E-04 & 2.46E-04 & 2.76E-04 & 2.20E-04 & 2.64E-04 & 2.55E-04 & 2.51E-04 & 1.58E-04 & 1.93E-04 & 2.95E-04 & --- \\
\hline
erru($t_k$)
& 8.08E-03 & 6.03E-03 & 5.90E-03 & 6.15E-03 & 6.41E-03 & 6.50E-03 & 5.88E-03 & 5.64E-03 & 6.48E-03 & 6.53E-03 & --- \\
errgu($t_k$)
& 2.41E-02 & 2.18E-02 & 2.13E-02 & 2.13E-02 & 2.14E-02 & 2.20E-02 & 2.07E-02 & 2.01E-02 & 1.98E-02 & 2.00E-02 & --- \\
errdivu($t_k$)
& 4.11E-02 & 2.49E-02 & 2.24E-02 & 2.38E-02 & 2.88E-02 & 3.09E-02 & 2.48E-02 & 2.30E-02 & 2.19E-02 & 2.15E-02 & --- \\
errp($t_k$)
& --- & --- & --- & --- & --- & --- & --- & --- & --- & --- & 1.93E-02 \\
        \hline
\end{tabular}}
		\caption{Error comparison.}
        \label{table:abc low beta}
\end{table}

\vspace{-0.5cm}

\subsection{Comparison with the deep Galerkin method (DGM)} 
\label{fjkd13}
 In this section, we compare the output of our method
 applied to the Taylor-Green vortex to that of
 the deep Galerkin method which has been developed in \cite{sirignano2018dgm}
 using neural networks.   
 In the following simulations we take $t=0$, $T=1/4$, and use the same number
 of neural network epochs as our deep branching (DB) algorithm,
 i.e. 20,000 epochs, and the computation times are comparable,
 as seen in Table~\ref{time_comparison}. 
 Note that the pre-computation of $p(T,x)$ is part of the 
 terminal boundary condition, and can be re-used for a different equation.

\begin{table}[H]
    \centering
            {\begin{tabular}{|r|c|c|} 
        \hline
         & Deep Branching (Taylor-Green) & DGM (Taylor-Green) \\ 
        \hline
        $p(T,x)$ & 1200s & \multirow{2}{*}{2400s} \\ 
        \cline{1-2} 
         $u(t,x)$ & 1300s &  \\ 
        \hline
\end{tabular}}
\caption{Comparison of computation times in seconds.}
\label{time_comparison} 
\end{table}

\vspace{-0.4cm}

\begin{figure}[H]
  \centering
 \begin{subfigure}[b]{0.45\textwidth}
    \includegraphics[width=\linewidth]{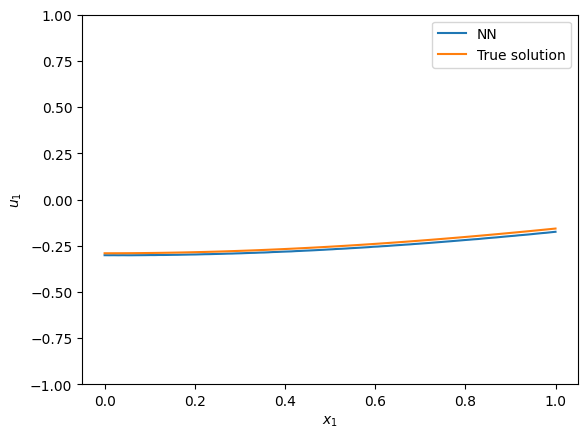}
    \caption{Comparison for $u_1(x_1,1)$.}
 \end{subfigure}
  \begin{subfigure}[b]{0.45\textwidth}
    \includegraphics[width=\linewidth]{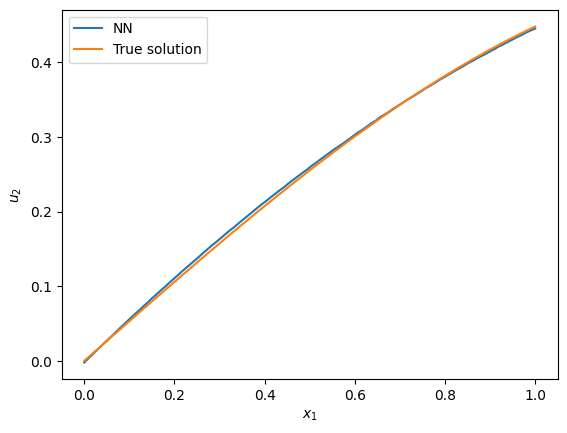} 
    \caption{Comparison for $u_2(x_1,1)$.}
 \end{subfigure}
   \vskip-0.1cm
   \caption{Comparison of DGM and \eqref{eq:taylor green} with space-time boundary condition and $x_2 = 1$.}
\label{fig:dgm1} 
\end{figure}

\noindent
In Figure~\ref{fig:dgm1}, we start with boundary conditions given by 
 \eqref{eq:taylor green} on the space-time domain $[0,1]^2 \times [0,T]$
 used in \cite{li-yue-zhang-duan}. 
 
\vspace{-0.4cm} 

\begin{figure}[H]
  \centering
 \begin{subfigure}[b]{0.45\textwidth}
    \includegraphics[width=\linewidth]{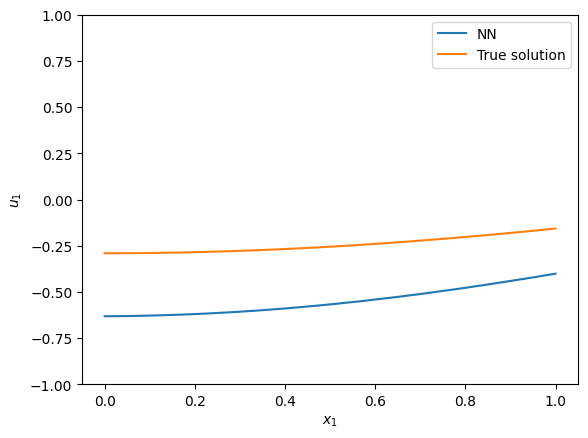}
    \caption{Comparison for $u_1(x_1,1)$.}
 \end{subfigure}
  \begin{subfigure}[b]{0.45\textwidth}
    \includegraphics[width=\linewidth]{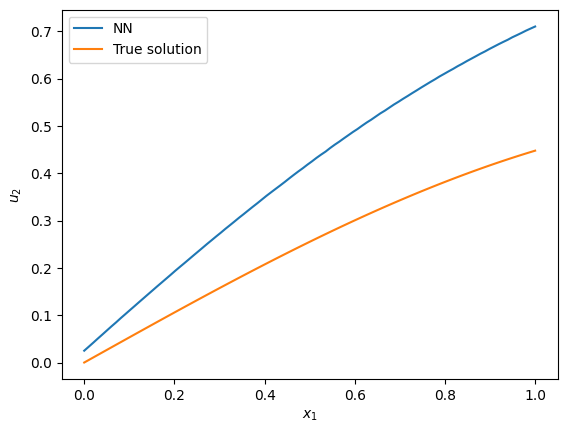} 
    \caption{Comparison for $u_2(x_1,1)$.}
 \end{subfigure}
   \vskip-0.1cm
   \caption{Comparison of DGM and \eqref{eq:taylor green} with terminal boundary condition and $x_2 = 1$.}
\label{fig:dgm2} 
\end{figure}

\noindent 
Next, in Figure~\ref{fig:dgm2} we only use a spatial boundary condition
on $[0,1]^2$ at the terminal time $T$, and we observe that accuracy
of the output is lost. 

\vskip-0.4cm

\begin{figure}[H]
  \centering
 \begin{subfigure}[b]{0.45\textwidth}
    \includegraphics[width=\linewidth]{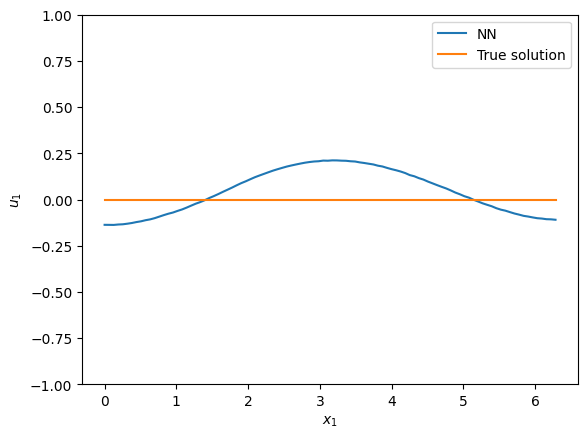}
    \caption{Comparison for $u_1(x_1,2\pi )$.}
 \end{subfigure}
  \begin{subfigure}[b]{0.45\textwidth}
    \includegraphics[width=\linewidth]{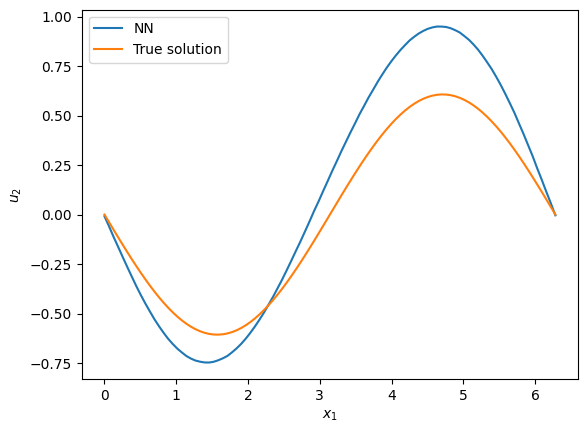} 
    \caption{Comparison for $u_2(x_1,2\pi )$.}
 \end{subfigure}
   \vskip-0.1cm
   \caption{Comparison of DGM and \eqref{eq:taylor green} with space-time boundary condition and $x_2 = 2\pi$.}
\label{fig:dgm3} 
\end{figure}

\begin{figure}[H]
  \centering
 \begin{subfigure}[b]{0.45\textwidth}
    \includegraphics[width=\linewidth]{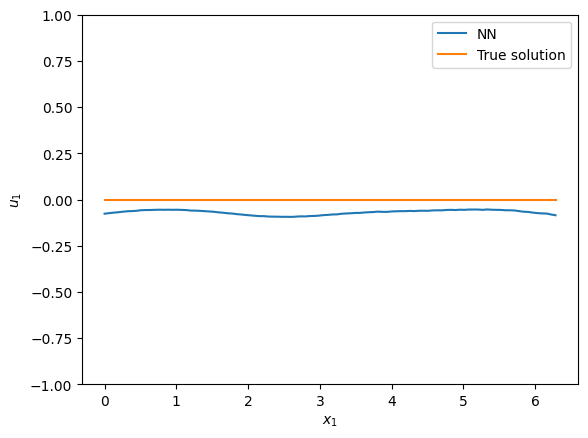}
    \caption{Comparison for $u_1(x_1,2\pi )$.}
 \end{subfigure}
  \begin{subfigure}[b]{0.45\textwidth}
    \includegraphics[width=\linewidth]{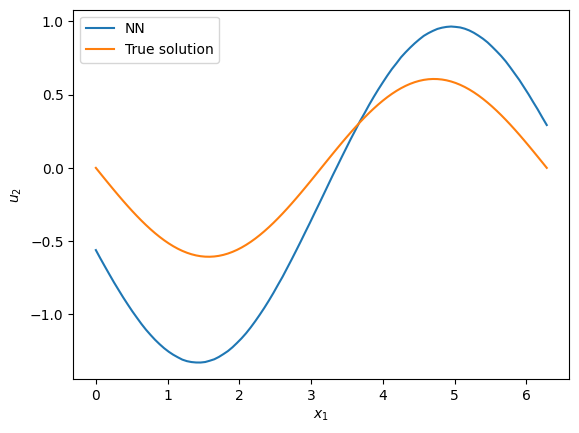} 
    \caption{Comparison for $u_2(x_1,2\pi )$.}
 \end{subfigure}
   \vskip-0.1cm
   \caption{Comparison of DGM and \eqref{eq:taylor green} with terminal boundary condition and $x_2 = 2\pi$.}
\label{fig:dgm4} 
\end{figure}

\noindent 
To conclude our assessment of the DGM method to the Taylor-Green vortex,
in Figures~\ref{fig:dgm3} and \ref{fig:dgm4} 
we extend the domain $[0,1]^2$ used in \cite{li-yue-zhang-duan}
to $[0,2\pi]^2$ as in Section~\ref{sub1}, and we observe that accuracy
is lost in this case, for both the space-time boundary condition on
 $[0,2\pi]^2\times [0,T]$ and the terminal boundary condition on
 $[0,2\pi]^2$ at time $T$. 

\vskip-0.4cm

\vskip-0.4cm

\subsection{Rotating flows} 
In this section we propose two other examples of dimensional flows on $\real^2$
that can be
solved with vanishing boundary conditions at infinity, by taking a
terminal condition $\phi$ of the form
$$
\left\{
\begin{array}{ll}
  \displaystyle
  \phi_1(x_1,x_2) = \frac{f'(x_2)}{f(x_2)} \exp \left( - \frac{g(x_1)}{f(x_2)} \right) 
  \medskip
  \\
  \displaystyle
  \phi_2(x_1,x_2) = \frac{g'(x_2)}{g(x_2)} \exp \left( - \frac{g(x_1)}{f(x_2)} \right), 
\end{array}
\right.
$$
which satisfies the divergence-free condition
${\rm div \ \! } \phi (x_1,x_2) = 0$, $(x_1,x_2)\in \real^2$.
This yields two-dimensional quiver velocity plots at different times
with $T=100$ in Figures~\ref{q1} and \ref{q2} below. 

\begin{figure}[H]
  \centering
 \begin{subfigure}[b]{0.45\textwidth}
    \includegraphics[width=\linewidth]{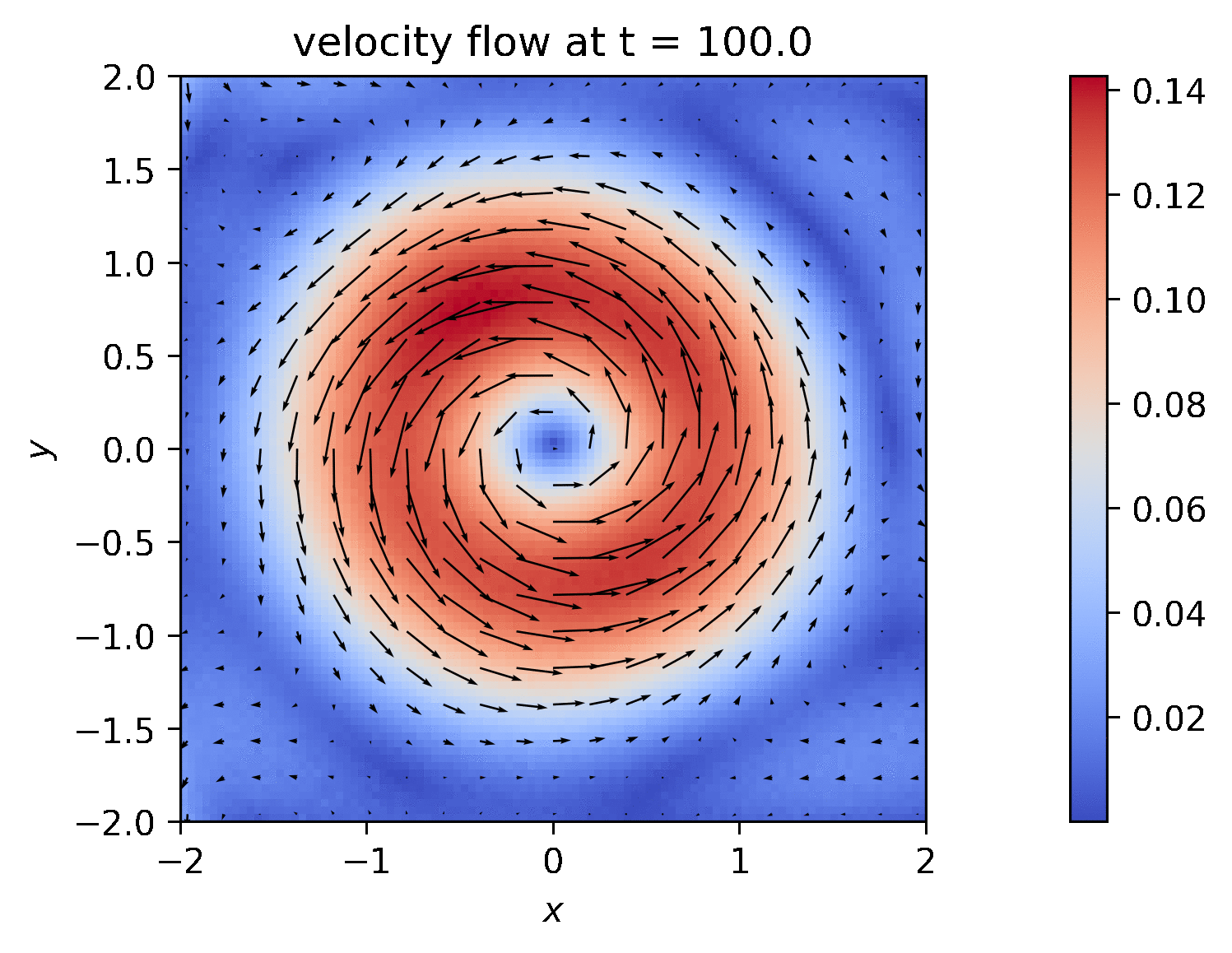} 
    \caption{$t=100$.}
 \end{subfigure}
  \begin{subfigure}[b]{0.45\textwidth}
    \includegraphics[width=\linewidth]{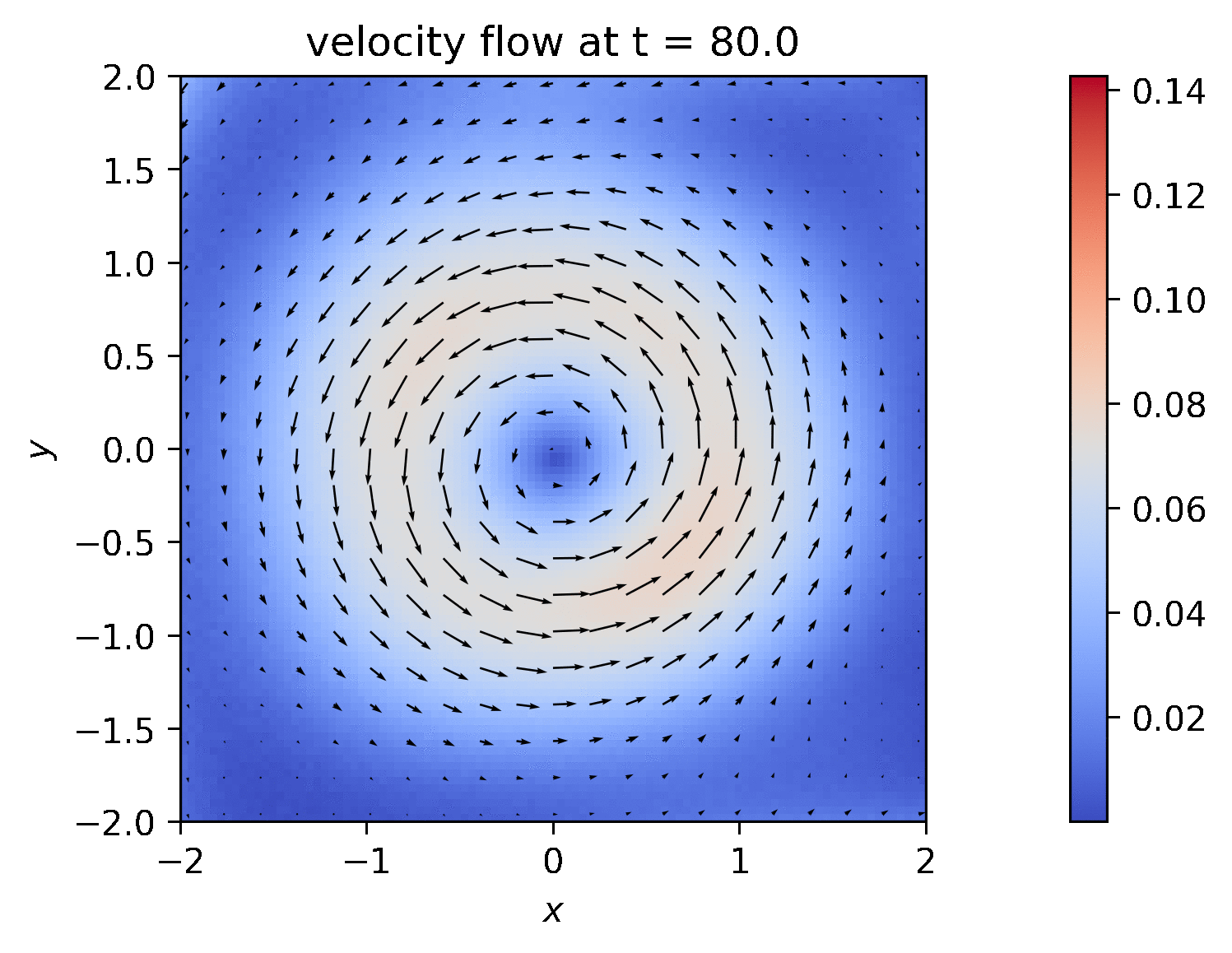} 
    \caption{$t=80$.}
 \end{subfigure}
   \vskip-0.1cm
   \caption{Case
  $f(x_2) = 1+x_2^2$, \ $g(x_1) = 1/(1+x_1^2)$.}      
\label{q1} 
\end{figure}

\vskip-0.4cm

\begin{figure}[H]
  \centering
 \begin{subfigure}[b]{0.45\textwidth}
    \includegraphics[width=\linewidth]{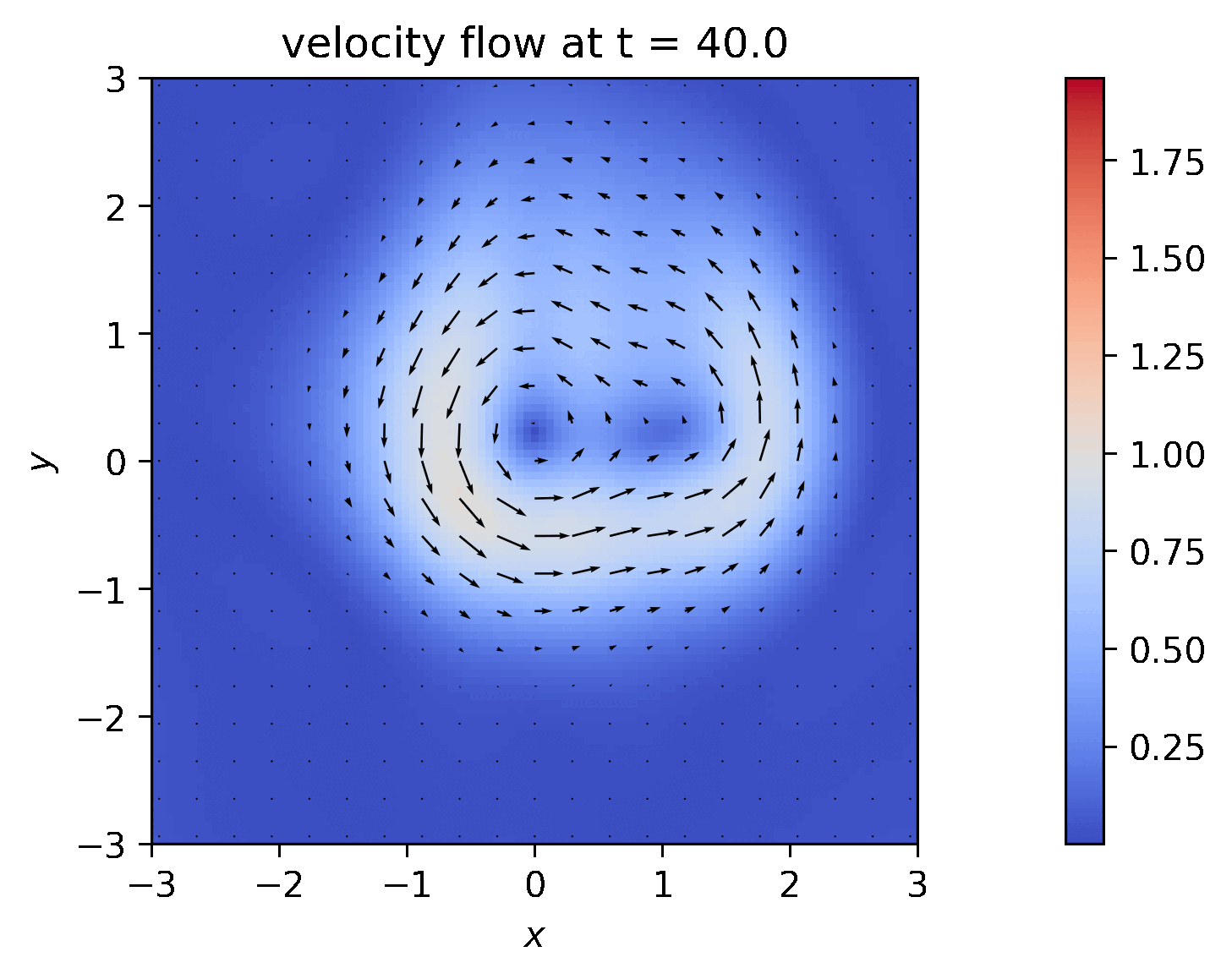} 
    \caption{$t=40$.}
 \end{subfigure}
  \begin{subfigure}[b]{0.45\textwidth}
    \includegraphics[width=\linewidth]{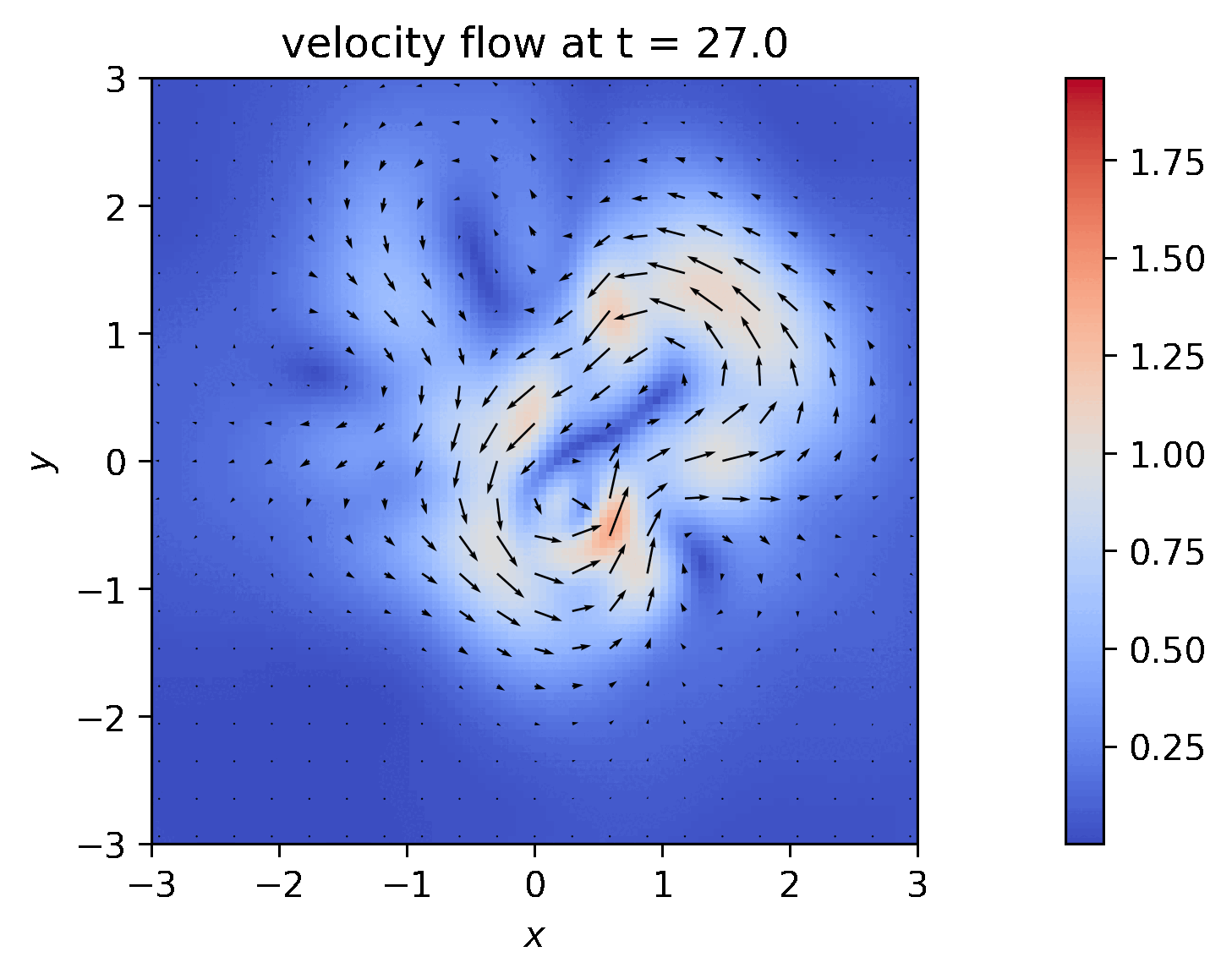}
    \caption{$t=27$.}
 \end{subfigure}
   \vskip-0.1cm
   \caption{Case
  $f(x_2) = (2+\sin (x_2) )/(1+x_2^2)$, \ $g(x_1) = \re^{x_1^2} / (2 + x_1^3 + x_1^4)$.}      
\label{q2} 
\end{figure}

\vskip-0.4cm

\appendix

\section{Branching solution of PDE systems}
In this section we present the extension of
the arguments of \cite{penent2022fully}, \cite{nguwipenentprivault}
to systems of partial differential equations,
 which leads to the probabilistic representation \eqref{eq:feynman kac}.
 The following proof uses the notation of Algorithm~\ref{alg:coding tree}.
 \begin{Proofy}\ref{p1}.
\noindent
$(i)$
     Consider $c\in \mathcal{C}$ a code of the form $c=(\partial_\lambda f)^*$.
From the Fa\`a di Bruno formula \eqref{eq:fdb} applied to the
function $f_{\beta_r}$, for $g \in {\cal C}^{\infty}(\real^d)$ we have
\begin{align}
\nonumber
& \partial_t g^*(u) + \nu \Delta g^*(u)
  \\
\nonumber
  & =   \sum\limits_{w = 1}^n
      \partial_{\bar{\alpha}^w}
    \left(\partial_t u_{\beta_w} + \nu \Delta u_{\beta_w} \right)
      (\partial_{\bm{1}_w} g)^*
  + \nu \sum\limits_{i = 1}^n \sum\limits_{j = 1}^n \sum\limits_{k = 1}^d
    \left(\partial_{\bar{\alpha}^i + \bm{1}_k} u_{\beta_i}\right)
    \left(\partial_{\bar{\alpha}^j + \bm{1}_k} u_{\beta_j}\right)
        (\partial_{\bm{1}_i + \bm{1}_j} g)^*
  \\
\nonumber
  &=
  \sum\limits_{w = 1}^{q}
      (\partial_{\bm{1}_w} g)^*
  \partial_{\bar{\alpha}^w} \left(\partial_t + \nu \Delta \right) u_0
    - \sum\limits_{w = q+1}^n
      (\partial_{\bm{1}_w} g)^*
    \partial_{\bar{\alpha}^w} ( f_{\beta_w}^*(u) )
  \\
\label{eq:partial t and Laplacian}
  &  \qquad \qquad \qquad
  + \nu \sum\limits_{i = 1}^n \sum\limits_{j = 1}^n \sum\limits_{k = 1}^d
    \left(\partial_{\bar{\alpha}^i + \bm{1}_k} u_{\beta_i}\right)
    \left(\partial_{\bar{\alpha}^j + \bm{1}_k} u_{\beta_j}\right)
        (\partial_{\bm{1}_i + \bm{1}_j} g)^*
  \\
\nonumber
  &=
  - \sum\limits_{w = q+1}^n
      (\partial_{\bm{1}_w} g)^*
    \left(\prod_{i = 1}^d \alpha^w_i! \right)
    \sum\limits_{\substack{1 \leq \lambda_1 + \cdots + \lambda_n \leq \abs{\bar{\alpha}^w} \\
                           1 \leq s \leq \abs{\bar{\alpha}^w}}}
    (\partial_{\lambda} f_{\beta_w})^*
    \sum\limits_{\substack{1 \leq \abs{k_1}, \dots, \abs{k_s}, \
                           0 \prec l^1 \prec \cdots \prec l^s \\
                           k^i_1 + \cdots + k^i_s = \lambda_i, \
                           i = 1, \dots, n \\
                           \abs{k_1}l_j^1 + \cdots + \abs{k_s}l_j^s = \alpha^w_j, \
                           j = 1, \dots, d
                           }}
    \prod_{\substack{1 \leq i \leq n \\
                     1 \leq r \leq s
                    }}
    \frac{(\partial_{l^r + \bar{\alpha}^i} u_{\beta_i})^{k_r^i}}{k_r^i!
                \left(l_1^r! \cdots l_d^r!\right)^{k_r^i}}
  \\
\nonumber
  &  \qquad \qquad \qquad
      + \sum\limits_{r = 1}^{q}
          (\partial_{\bm{1}_r} g)^*
      \partial_{\bar{\alpha}^r} \left(\partial_t + \nu \Delta \right) u_0
      + \nu \sum\limits_{i = 1}^n \sum\limits_{j = 1}^n \sum\limits_{k = 1}^d
    \left(\partial_{\bar{\alpha}^i + \bm{1}_k} u_{\beta_i}\right)
    \left(\partial_{\bar{\alpha}^j + \bm{1}_k} u_{\beta_j}\right)
        (\partial_{\bm{1}_i + \bm{1}_j} g)^*.
\end{align}
 Rewriting the above equation in integral form yields
 \begin{align}
   \label{fjkld133} 
   & g^*(u)(t,x) = \int_{\real^d} \varphi_{2\nu} (T-t,y-x) g^* (\phi)(y) dy
   \\
   \nonumber
   &
   + \int_t^T \int_{\real^d} \varphi_{2\nu} (s-t,y-x)
   \Biggl(
      - \nu \sum\limits_{i = 1}^n \sum\limits_{j = 1}^n \sum\limits_{k = 1}^d
     \big(\partial_{\bar{\alpha}^i + \bm{1}_k} u_{\beta_i}(s, y)\big)
     \big(\partial_{\bar{\alpha}^j + \bm{1}_k} u_{\beta_j}(s, y)\big)
        (\partial_{\bm{1}_i + \bm{1}_j} g)^*
    \\
\nonumber
      &
      - \sum\limits_{r = 1}^{q}
          (\partial_{\bm{1}_r} g)^*
          \partial_{\alpha^r} \left(\partial_t + \nu \Delta \right) u_0(s, y)
    \\
\nonumber
      &
   + \sum\limits_{w = q+1}^n
      (\partial_{\bm{1}_w} g)^*
    \left(\prod_{i = 1}^d {\alpha}^w_i! \right)
    \sum\limits_{\substack{1 \leq \lambda_1 + \cdots + \lambda_n \leq \abs{\bar{\alpha}^w} \\
                           1 \leq s \leq \abs{\bar{\alpha}^w}}}
    \hskip-0.3cm
    (\partial_{\lambda} f_{\beta_w})^*
        \hskip-0.9cm
    \sum\limits_{\substack{1 \leq \abs{k_1}, \dots, \abs{k_s}, \
                           0 \prec l^1 \prec \cdots \prec l^s \\
                           k^i_1 + \cdots + k^i_s = \lambda_i, \
                           i = 1, \dots, n \\
                           \abs{k_1}l_j^1 + \cdots + \abs{k_s}l_j^s = \bar{\alpha}^w_j, \
                           j = 1, \dots, d
                           }}
    \prod_{\substack{1 \leq i \leq n \\
                     1 \leq r \leq s
                    }}
 \frac{\left(\partial_{l^r + \bar{\alpha}^i} u_{\beta_i}(s, y)\right)^{k_r^i}}{k_r^i!
                \left(l_1^r! \cdots l_d^r!\right)^{k_r^i}}
   \Bigg) dy ds,
\end{align}
 which shows the equation
\begin{equation}
\label{s1-}
 c(u)(t,x) = \int_{\real^d} \varphi_{2\nu} (T-t,y-x) c ( u)(T,y) dy
+
\sum_{Z \in \mathcal{M}(c)}
\int_t^T \int_{\real^d} \varphi_{2\nu} (s-t,y-x)
\prod_{z \in Z}  z(u)(s,y) dy ds,
\end{equation}
$(t,x)\in [0,T]\times \real$,
for any code $c\in \mathcal{C}$ of the form $c=(\partial_\lambda f)^*$.
Also, \eqref{s1-} holds directly 
from \eqref{eq:integral u0} for the code $c = {\rm Id}_i$,
$i=1,\ldots , d$. 

\medskip

\noindent
$(ii)$
For $c = ( \partial_{\mu}, 0)$,
by \eqref{eq:integral u0}
we have
\begin{eqnarray}
  \lefteqn{
    \partial_\mu u_0(t, x)
  =
  \nonumber
    \frac{\Gamma (d/2)}{2 \pi^{d/2}}
    \int_{\real^d}
 \frac{ N(y)}{  \abs{y}^d}
    \partial_\mu
    f_0\big(
          \partial_{\bar{\alpha}^{q+1}}u_{\beta_{q+1}}(t,x+y)
        ,
        \ldots ,
                   \partial_{\bar{\alpha}^n}u_{\beta_n}(t,x+y)
                   \big)
                     dy
  }
  \\
  \nonumber
  & = &
   \int_{\real^d}
  N(y)
  \int_0^\infty
  \frac{(2 \pi s)^{-d/2}}{2s}
  \re^{-\abs{y}^2/(2s)}
  ds
  \partial_\mu
       f_0\big(
          \partial_{\bar{\alpha}^{q+1}}u_{\beta_{q+1}}(t,x+y)
        ,
        \ldots ,
                   \partial_{\bar{\alpha}^n}u_{\beta_n}(t,x+y)
                   \big)
  dy
  \\
\nonumber
  &= &
  \int_0^\infty
  \int_{\real^d}
  \varphi_1 (s, y)
  \frac{N(y)}{2s}
  \partial_\mu
           f_0\big(
          \partial_{\bar{\alpha}^{q+1}}u_{\beta_{q+1}}(t,x+y)
        ,
        \ldots ,
                   \partial_{\bar{\alpha}^n}u_{\beta_n}(t,x+y)
                   \big)
  dy
  ds
  \\
  \label{eq:poisson integration calculation}
  &= &
  \sum\limits_{Z \in \mathcal{M}((\partial_\mu, 0))}
  \int_0^\infty
  \int_{\real^d}
  \varphi_1 (s, y)
  \frac{N(y)}{2s}
  \prod\limits_{z \in Z}
  z(u)(t, x + y)
  dy
  ds, 
\end{eqnarray}
 which shows that 
\begin{equation}
  \label{fjkldf21-0}
  c(u)(t,x) =
  \sum\limits_{Z \in \mathcal{M}(c)}
  \int_0^\infty
  \int_{\real^d}
  \varphi_1 (s, y)
  \frac{N(y)}{2s}
  \prod\limits_{z \in Z}
  z(u)(t, x + y)
  dy
  ds,
\end{equation}
for the code $c = ( \partial_{\mu}, 0)$.
Also, \eqref{fjkldf21-0} holds directly for the code $c = {\rm Id}_0$
from \eqref{eq:integral u0}. 

\medskip

\noindent
$(iii)$ Next, for the code $c = (\partial_{\mu}, -1)$, from
 \eqref{eq:partial t and Laplacian} applied to $g=f_0$,
 we have
\begin{align}
  \nonumber
  &
  \partial_\mu (\partial_t + \nu \Delta) u_0(t, x)
    \\
  \nonumber
  &=
  \int_0^\infty
  \int_{\real^d}
  \varphi_1 (s, y)
  \frac{N(y)}{2s}
    \partial_\mu
  (\partial_t + \nu \Delta)
              f_0\big(
          \partial_{\bar{\alpha}^{q+1}}u_{\beta_{q+1}}(t,x+y)
        ,
        \ldots ,
                   \partial_{\bar{\alpha}^n}u_{\beta_n}(t,x+y)
                   \big)
  dy
  ds
  \\
  \nonumber
  &=
  \int_0^\infty
  \int_{\real^d}
  \varphi_1 (s, y)
  \frac{N(y)}{2s}
  \\
\nonumber
    &  \qquad
  \partial_\mu
  \left(
  \nu \sum\limits_{i = q+1}^n \sum\limits_{j = q+1}^n \sum\limits_{k = 1}^d
    \left(\partial_{\bar{\alpha}^i + \bm{1}_k} u_{\beta_i}\right)
    \left(\partial_{\bar{\alpha}^j + \bm{1}_k} u_{\beta_j}\right)
        (\partial_{\bm{1}_i + \bm{1}_j} f_0)^*
  - \sum\limits_{i = q+1}^n
      (\partial_{\bm{1}_i} f_0)^*
    \partial_{\bar{\alpha}^i} f_{\beta_i}^*(u)
  \right)
  dy
  ds
  \\
  \nonumber
  &=
  \int_0^\infty
  \int_{\real^d}
  \varphi_1 (s, y)
  \frac{N(y)}{2s}
  \\
\nonumber
    &  \qquad
\left(
\sum\limits_{\substack{0 \leq \gamma_i \leq \ell_i \leq \mu_i \\
                         1 \leq i \leq d}}
  \left(
      \nu
      \prod\limits_{r = 1}^d
      {\mu_r \choose \ell_r}
      {\ell_r \choose \gamma_r}
  \right)
  \sum\limits_{i = q+1}^n \sum\limits_{j = q+1}^n \sum\limits_{k = 1}^d
    \left(\partial_{\mu - \ell + \bar{\alpha}^i + \bm{1}_k} u_{\beta_i}\right)
    \left(\partial_{\ell - \gamma + \bar{\alpha}^j + \bm{1}_k} u_{\beta_j}\right)
        \partial_{\gamma}(\partial_{\bm{1}_i + \bm{1}_j} f_0)^*
        \right.
        \\
\nonumber
    &  \qquad \qquad
\left.
- \sum\limits_{\substack{0 \leq \ell_i \leq \mu_i \\
                         1 \leq i \leq d}}
  \left(
      \prod\limits_{r = 1}^d
      {\mu_r \choose \ell_r}
  \right)
  \sum\limits_{i = q+1}^n
  \partial_{\ell} (\partial_{\bm{1}_i} f_0)^*
    \partial_{\mu - \ell + \bar{\alpha}^i} f_{\beta_i}^*(u)
    \right)
    dy
  ds
  \\
  \nonumber
  &=
  \int_0^\infty
  \int_{\real^d}
  \varphi_1 (s, y)
  \frac{N(y)}{2s}
  \\
  \nonumber
  &  
  \left(
  \sum\limits_{\substack{0 \leq \gamma_i \leq \ell_i \leq \mu_i \\
                         1 \leq i \leq d}}
  \left(
      \nu
      \prod\limits_{r = 1}^d
      {\mu_r \choose \ell_r}
      {\ell_r \choose \gamma_r}
  \right)
  \sum\limits_{i = q+1}^n \sum\limits_{j = q+1}^n \sum\limits_{k = 1}^d
    \left(\partial_{\mu - \ell + \bar{\alpha}^i + \bm{1}_k} u_{\beta_i}\right)
    \left(\partial_{\ell - \gamma + \bar{\alpha}^j + \bm{1}_k} u_{\beta_j}\right)
        \partial_{\gamma}(\partial_{\bm{1}_i + \bm{1}_j} f_0)^*
        \right.
        \\
  \nonumber
  & 
  - \sum\limits_{\substack{0 \leq \ell_i \leq \mu_i \\
                         1 \leq i \leq d}}
  \left(
      \prod\limits_{w = 1}^d
      {\mu_w \choose \ell_w}
      \ell_w!
  \right)
\hskip-0.6cm
  \sum\limits_{q < i \leq n
    \atop
        {1 \leq \lambda_1 + \cdots + \lambda_n \leq \abs{\ell}
          \atop
          1 \leq s \leq \abs{\ell}
        }
        }
    \hskip-0.8cm
  \partial_{\mu - \ell + \bar{\alpha}^i} f_{\beta_i}^*(u)
   (\partial_{\lambda + \bm{1}_i} f_0)^* (u)
   \hskip-1.1cm
   \sum\limits_{\substack{1 \leq \abs{k_1}, \dots, \abs{k_s}, \
                           0 \prec l^1 \prec \cdots \prec l^s \\
                           k^i_1 + \cdots + k^i_s = \lambda_i, \
                           i = 1, \dots, n \\
                           \abs{k_1}l_j^1 + \cdots + \abs{k_s}l_j^s = \ell_j, \
                           j = 1, \dots, d
                           }}
   \hskip-0.1cm
   \left.
  \prod_{\substack{1 \leq i \leq n \\
                     1 \leq r \leq s
                    }}
                    \frac{(\partial_{l^r + \bar{\alpha}^i} u_{\beta_i} )^{k_r^i}}{k_r^i!
                \left(l_1^r! \cdots l_d^r!\right)^{k_r^i}}
                  \right)
                dy
  ds
  \\
 \label{jkldf09} 
  \\
  \nonumber
  \\
  \label{eq:poisson integration calculation-2}
  &=
  \sum\limits_{Z \in \mathcal{M}((\partial_\mu, -1))}
  \int_0^\infty
  \int_{\real^d}
  \varphi_1 (s, y)
  \frac{N(y)}{2s}
  \prod\limits_{z \in Z}
  z(u)(t, x + y)
  dy
  ds, 
\end{align}
according to the definition of $\mathcal{M}((\partial_\mu, -1))$. 
Hence, we have shown that
\begin{equation}
  \label{fjkldf21}
  c(u)(t,x) =
  \sum\limits_{Z \in \mathcal{M}(c)}
  \int_0^\infty
  \int_{\real^d}
  \varphi_1 (s, y)
  \frac{N(y)}{2s}
  \prod\limits_{z \in Z}
  z(u)(t, x + y)
  dy
  ds
\end{equation}
for the code $c = (\partial_{\mu}, -1)$.
 
\noindent
$(iv)$
 By the Fa\`a di Bruno formula \eqref{eq:fdb}, 
 Equation~\eqref{s1-} is also satisfied by
 $c = ( \partial_{\mu}, i)$ for $i=1,\ldots , d$,
 since ${\cal M}\left((\partial_{\mu}, i)\right)
 = {\rm fdb} (\mu, f_i, \emptyset )$.

 \medskip
  
 \noindent
 For any $c\in \mathcal{C}$, we now let
$$
   u_c(t,x) := \E [ \mathcal{H}(t,x,c) ],
   \qquad
   (t,x) \in [0,T] \times \real^d.
$$
\noindent
$(v)$
 Starting from a code of the form
  $c=(\partial_\lambda f)^*$ or
 $c = ( \partial_{\mu}, i)$ for $i=1,\ldots , d$,
 we draw a sample of $I_c$ uniformly in $\mathcal{M}(c)$.
As each code in the tuple $I_c$ yields a new branch
at time $\tau$, we have 
\begin{align}
\nonumber
&
 u_c (t,x) =  \E \big[ \mathcal{H}(t,x,c)\mathbbm{1}_{\{ t + \tau > T  \}} + \mathcal{H}(t,x,c) \mathbbm{1}_{\{ t + \tau \leq T \}} \big]
\\
\nonumber
&= \E \left[ \frac{c(u)(T, x + W_{2\nu (T-t)}) }{\widebar{F}(T-t)} \mathbbm{1}_{\{ t + \tau > T \}}
+
  \mathbbm{1}_{\{ t + \tau \leq T \}}
  \sum_{Z\in {\cal M}(c)}
  {\bf 1}_{\{ I_c = Z\}}
r_c \frac{ \prod_{z\in Z} u_z ( t + \tau  , x + W_{2\nu \tau} )}{ \rho( \tau )}
\right]
\\
\label{fjhkdsf}
  & =
  \int_{-\infty}^\infty \varphi_{2\nu} (T-t,y-x) c(u)(T,y) dy
+
\sum_{Z \in \mathcal{M}(c)}
\int_t^T \int_{-\infty}^\infty \varphi_{2\nu} (s-t,y-x)
\prod_{z \in Z}  u_z (s,y) dy ds,
\end{align}
which yields the same system of equations as \eqref{s1-}.

\noindent
$(vi)$
Similarly, starting 
 a code of the form $c = ( \partial_{\mu}, 0)$ or $c = (\partial_{\mu}, -1)$
we draw a sample of $I_c$ uniformly in $\mathcal{M}(c)$ with probability $1/r_c$, where
 $r_c$ is the size of $\mathcal{M}(c)$.
As each code in the tuple $I_c$ yields a new branch at time $\widetilde{\tau}$, we obtain
\begin{align}
\nonumber
 u_c (t,x) & = \E [ \mathcal{H}(t,x,c) ]
\\
\nonumber
&= \E \left[
  \frac{r_c
 N(W_{\widetilde{\tau}})
  }{ 2 \widetilde{\tau} \widetilde{\rho}(\widetilde{\tau})}
  \sum_{Z\in {\cal M}(c)}
  {\bf 1}_{\{ I_c = Z\}}
  \prod_{z\in Z} u_z (t, x + W_{\widetilde{\tau}})
\right]
\\
\label{fjhkdsf2}
    & =
  \sum\limits_{Z \in \mathcal{M}(c)}
  \int_0^\infty
  \int_{\real^d}
  \varphi_1 (s, y)
  \frac{N(y)}{2s}
  \prod\limits_{z \in Z}
  z(u)(t, x + y)
  dy
  ds,
\end{align}
which coincides with \eqref{fjkldf21-0} or \eqref{fjkldf21},
respectively for 
$c = ( \partial_{\mu}, 0)$ and $c = (\partial_{\mu}, -1)$.

\noindent
$(vii)$
From \eqref{fjhkdsf}-\eqref{fjhkdsf2} and
\eqref{s1-}-\eqref{fjkldf21-0}-\eqref{fjkldf21}
we conclude that for any code $c\in \mathcal{C}$,
$u_c(t,x)$ and $c(u)(t,x)$ satisfy the same system of equations
\eqref{s1--}. 
As by assumption the system \eqref{s1--} has a unique solution
we conclude that $(c (u))_{c\in \mathcal{C} } = (u_c)_{c\in \mathcal{C} }$, and
therefore
$$
u_c (t,x) = \E [ \mathcal{H}(t,x,c) ]
=
c(u)(t,x), \qquad
(t,x) \in [0,T] \times \real,
\quad c \in \mathcal{C}.
$$
In particular, for $c={\rm Id}_i$ this yields 
$$
u_i (t,x) =
{\rm Id}_i(u)(t,x)
= \E [ \mathcal{H}(t,x,{\rm Id}_i) ],  \qquad
(t,x) \in [0,T] \times \real,
\quad i=0,1,\ldots , d, 
$$
which is \eqref{eq:feynman kac}. 
\end{Proofy}

\footnotesize

\newcommand{\etalchar}[1]{$^{#1}$}
\def\cprime{$'$} \def\polhk#1{\setbox0=\hbox{#1}{\ooalign{\hidewidth
  \lower1.5ex\hbox{`}\hidewidth\crcr\unhbox0}}}
  \def\polhk#1{\setbox0=\hbox{#1}{\ooalign{\hidewidth
  \lower1.5ex\hbox{`}\hidewidth\crcr\unhbox0}}} \def\cprime{$'$}


\begin{thebibliography}{HLOT{\etalchar{+}}19}

\bibitem[AB10]{belopolskaya}
S.~Albeverio and Ya. Belopolskaya.
\newblock Generalized solutions of the {C}auchy problem for the
  {N}avier-{S}tokes system and diffusion processes.
\newblock {\em Cubo}, 12(2):77--96, 2010.

\bibitem[APFC17]{angeli}
P.-E. Angeli, M.-A. Puscas, G.~Fauchet, and A.~Cartalade.
\newblock {FVCA8} {B}enchmark for the {S}tokes and {N}avier–{S}tokes
  equations with the {T}rio{CFD} code-benchmark session.
\newblock In {\em FVCA 2017: Finite Volumes for Complex Applications VIII -
  Methods and Theoretical Aspects}, volume 199 of {\em Springer Proceedings in
  Mathematics \& Statistics}, pages 181--202. Springer Verlag, 2017.

\bibitem[Arn65]{arnold1965topologie}
V.~Arnol{\cprime d}.
\newblock Sur la topologie des \'{e}coulements stationnaires des fluides
  parfaits.
\newblock {\em C. R. Acad. Sci. Paris}, 261:17--20, 1965.

\bibitem[Bor17]{borodin}
A.N. Borodin.
\newblock {\em Stochastic processes}.
\newblock Probability and its Applications. Birkh\"{a}user/Springer, Cham,
  2017.
\newblock Original Russian edition published by LAN Publishing, St. Petersburg,
  2013.

\bibitem[CC07]{cipriano}
F.~Cipriano and A.B. Cruzeiro.
\newblock {N}avier-{S}tokes equation and diffusions on the group of
  homeomorphisms of the torus.
\newblock {\em Comm. Math. Phys.}, 275:255--269, 2007.

\bibitem[Chi70]{childress1970new}
S.~Childress.
\newblock New solutions of the kinematic dynamo problem.
\newblock {\em J. Math. Phys.}, 11(10):3063--3076, 1970.

\bibitem[CS96]{constantine}
G.M. Constantine and T.H. Savits.
\newblock A multivariate {F}aa di {Bruno} formula with applications.
\newblock {\em Trans. Amer. Math. Soc.}, 348(2):503--520, 1996.

\bibitem[CS09]{cruzeiroshamarova}
A.B. Cruzeiro and E.~Shamarova.
\newblock Navier-{S}tokes equations and forward-backward {SDE}s on the group of
  diffeomorphisms of a torus.
\newblock {\em Stochastic Process. Appl.}, 119(12):4034--4060, 2009.

\bibitem[CSTV07]{touzi}
P.~Cheridito, H.M. Soner, N.~Touzi, and N.~Victoir.
\newblock Second-order backward stochastic differential equations and fully
  nonlinear parabolic {PDE}s.
\newblock {\em Comm. Pure Appl. Math.}, 60(7):1081--1110, 2007.

\bibitem[DQT15]{delbaen2015forward}
F.~Delbaen, J.~Qiu, and S.~Tang.
\newblock Forward-backward stochastic differential systems associated to
  {N}avier-{S}tokes equations in the whole space.
\newblock {\em Stochastic Process. Appl.}, 125(7):2516--2561, 2015.

\bibitem[FTW11]{fahim}
A.~Fahim, N.~Touzi, and X.~Warin.
\newblock A probabilistic numerical method for fully nonlinear parabolic
  {PDE}s.
\newblock {\em Ann. Appl. Probab.}, 21(4):1322--1364, 2011.

\bibitem[GZZ15]{guowenjie}
W.~Guo, J.~Zhang, and J.~Zhuo.
\newblock A monotone scheme for high-dimensional fully nonlinear {PDE}s.
\newblock {\em Ann. Appl. Probab.}, 25(3):1540--1580, 2015.

\bibitem[HJE18]{han2018solving}
J.~Han, A.~Jentzen, and W.~E.
\newblock Solving high-dimensional partial differential equations using deep
  learning.
\newblock {\em Proceedings of the National Academy of Sciences},
  115(34):8505--8510, 2018.

\bibitem[HL12]{henry-labordere2012}
P.~Henry-Labord\`ere.
\newblock Counterparty risk valuation: a marked branching diffusion approach.
\newblock Preprint arXiv:1203.2369, 2012.

\bibitem[HLOT{\etalchar{+}}19]{labordere}
P.~Henry-Labord\`ere, N.~Oudjane, X.~Tan, N.~Touzi, and X.~Warin.
\newblock Branching diffusion representation of semilinear {PDE}s and {M}onte
  {C}arlo approximation.
\newblock {\em Ann. Inst. H. Poincar\'e Probab. Statist.}, 55(1):184--210,
  2019.

\bibitem[HLZ20]{huangshuo}
S.~Huang, G.~Liang, and T.~Zariphopoulou.
\newblock An approximation scheme for semilinear parabolic {PDE}s with convex
  and coercive {H}amiltonians.
\newblock {\em SIAM J. Control Optim.}, 58(1):165--191, 2020.

\bibitem[Hor91]{hornik1991approximation}
K.~Hornik.
\newblock Approximation capabilities of multilayer feedforward networks.
\newblock {\em Neural networks}, 4(2):251--257, 1991.

\bibitem[HZRS16]{He16}
K.~He, X.~Zhang, S.~Ren, and J.~Sun.
\newblock Deep residual learning for image recognition.
\newblock In {\em Proceedings of the IEEE conference on computer vision and
  pattern recognition}, pages 770--778, 2016.

\bibitem[INW69]{inw}
N.~Ikeda, M.~Nagasawa, and S.~Watanabe.
\newblock Branching {M}arkov processes {I}, {II}, {III}.
\newblock {\em J. Math. Kyoto Univ.}, 8-9:233--278, 365--410, 95--160,
  1968-1969.

\bibitem[IS15]{ioffe2015batch}
S.~Ioffe and Ch. Szegedy.
\newblock Batch normalization: {Accelerating} deep network training by reducing
  internal covariate shift.
\newblock {\em Preprint arXiv:1502.03167}, 2015.

\bibitem[KB14]{kingma2014adam}
D.P. Kingma and J.~Ba.
\newblock Adam: {A} method for stochastic optimization.
\newblock {\em Preprint arXiv:1412.6980}, 2014.

\bibitem[LG20]{lejay2020forward}
A.~Lejay and H.M. Gonz{\'a}lez.
\newblock A forward-backward probabilistic algorithm for the incompressible
  {N}avier-{S}tokes equations.
\newblock {\em Journal of Computational Physics}, 420:109689, 19, 2020.

\bibitem[LYZD22]{li-yue-zhang-duan}
J.~Li, J.~Yue, W.~Zhang, and W.~Duan.
\newblock The deep learning {G}alerkin method for the general {S}tokes
  equations.
\newblock {\em J. Sci. Comput.}, 93(1):Paper No. 5, 20, 2022.

\bibitem[Mat21]{mmatsumoto}
M.~Matsumoto.
\newblock Application of {D}eep {G}alerkin {M}ethod to solve compressible
  {N}avier-{S}tokes equations.
\newblock {\em Trans. Japan Soc. Aero. Space Sci.}, 64(6):348--357, 2021.

\bibitem[MB02]{majda2002vorticity}
A.J. Majda and A.L. Bertozzi.
\newblock {\em Vorticity and incompressible flow}, volume~27 of {\em Cambridge
  Texts in Applied Mathematics}.
\newblock Cambridge University Press, Cambridge, 2002.

\bibitem[McK75]{hpmckean}
H.P. McKean.
\newblock Application of {B}rownian motion to the equation of
  {K}olmogorov-{P}etrovskii-{P}iskunov.
\newblock {\em Comm. Pure Appl. Math.}, 28(3):323--331, 1975.

\bibitem[NPP22a]{nguwipenentprivault}
J.Y. Nguwi, G.~Penent, and N.~Privault.
\newblock A deep branching solver for fully nonlinear partial differential
  equations.
\newblock Preprint arXiv:2203.03234, 17 pages, 2022.

\bibitem[NPP22b]{penent2022fully}
J.Y. Nguwi, G.~Penent, and N.~Privault.
\newblock A fully nonlinear {F}eynman-{K}ac formula with derivatives of
  arbitrary orders.
\newblock Preprint arXiv:2201.03882v3, 30 pages, 2022.

\bibitem[PP92]{pardouxpeng}
{\'E}.~Pardoux and S.~Peng.
\newblock Backward stochastic differential equations and quasilinear parabolic
  partial differential equations.
\newblock In {\em Stochastic partial differential equations and their
  applications ({C}harlotte, {NC}, 1991)}, volume 176 of {\em Lecture Notes in
  Control and Inform. Sci.}, pages 200--217. Springer, Berlin, 1992.

\bibitem[PP22]{penent4}
G.~Penent and N.~Privault.
\newblock Numerical evaluation of {ODE} solutions by {M}onte {C}arlo
  enumeration of {B}utcher series.
\newblock Preprint arXiv:2201.05998, to appear in BIT Numerical Mathematics,
  2022.

\bibitem[RPK19]{Rassi19}
M.~Raissi, P.~Perdikaris, and G.~E. Karniadakis.
\newblock Physics-informed neural networks: A deep learning framework for
  solving forward and inverse problems involving nonlinear partial differential
  equations.
\newblock {\em Journal of Computational Physics}, 378:686--707, 2019.

\bibitem[Sko64]{skorohodbranching}
A.V. Skorokhod.
\newblock Branching diffusion processes.
\newblock {\em Teor. Verojatnost. i. Primenen.}, 9:492--497, 1964.

\bibitem[SS18]{sirignano2018dgm}
J.~Sirignano and K.~Spiliopoulos.
\newblock {DGM}: {A} deep learning algorithm for solving partial differential
  equations.
\newblock {\em Journal of Computational Physics}, 375:1339--1364, 2018.

\bibitem[STZ12]{soner}
H.M. Soner, N.~Touzi, and J.~Zhang.
\newblock Wellposedness of second order backward {SDE}s.
\newblock {\em Probab. Theory Related Fields}, 153(1-2):149--190, 2012.

\bibitem[Tan13]{tanxiaolu}
X.~Tan.
\newblock A splitting method for fully nonlinear degenerate parabolic {PDE}s.
\newblock {\em Electron. J. Probab.}, 18:no. 15, 24, 2013.

\bibitem[TG37]{taylor1937mechanism}
G.I. Taylor and A.E. Green.
\newblock Mechanism of the production of small eddies from large ones.
\newblock {\em Proc. R. Soc. Lond. Ser. A Math. Phys. Eng. Sci.},
  158(895):499--521, 1937.

\end{thebibliography}
\end{document}